\documentclass{svmult}

\usepackage{makeidx}         
\usepackage{graphicx}        
\usepackage{multicol}        
\usepackage[active]{srcltx}
\usepackage{amsmath,amssymb}

\makeindex             

\newcommand{\bydef}{\stackrel{\mbox{\scriptsize def}}{=}}
\def\H{{\mathcal{H}}}
\def\K{{\mathcal{K}}}

\begin{document}

\title*{Inverse Problems for Semiconductors: Models and Methods}

\author{A.\,Leit\~ao\inst{1}    \and
        P.A.\,Markowich\inst{2} \and
        J.P.\,Zubelli\inst{3} }

\institute{
\mbox{}\!\!\!\!
Department of Mathematics, Federal University of St.\,Catarina,
P.O. Box 476, 88040-900 Florianopolis, Brazil
\texttt{\{aleitao@mtm.ufsc.br\}}
\and
\mbox{}\!\!\!\!
Department of Mathematics, University of Vienna, Boltzmanngasse 9,
A-1090 Vienna, Austria \texttt{\{peter.markowich@univie.ac.at\}}
\and
\mbox{}\!\!\!\!
IMPA, Estr.\,D.\,Castorina\,110, 22460-320\,Rio\,de\,Janeiro, Brazil%
\,\texttt{\{zubelli@impa.br\}} }
%
%
\maketitle
\setcounter{footnote}{0}

\begin{abstract}
We consider the problem of identifying discontinuous doping
profiles in semiconductor devices from data obtained by different
models connected to the voltage-current map. 
Stationary as well as transient settings are discussed and a framework
for the corresponding inverse problems is established.
Numerical implementations for the so-called {\em stationary unipolar and
stationary bipolar cases} show the effectiveness of a level set
approach to tackle the inverse problem.
\end{abstract}

\section{Introduction} \label{sec:1}

\paragraph{The mathematical model}

The starting point of the mathematical model discussed in this paper is
the system of {\em drift diffusion equations} (see (\ref{eq:VanR1}) --
(\ref{eq:VanR6}) below).
This system of equations, derived more than fifty years ago \cite{vRo50},
is the most widely used to describe semiconductor devices.
For the current state of technology, this system represents an accurate
compromise between efficient numerical solvability of the mathematical
model and realistic description of the underlying physics
\cite{Mar86,MRS90,Sel84}.

The name {\em drift diffusion equations} of semiconductors originates
from the type of dependence of the current densities on the carrier
densities and the electric field. The current densities are the sums of
drift terms and diffusion terms.
It is worth mentioning that, with the increased miniaturization of
semiconductor devices, one comes closer and closer to the limits of
validity of the drift diffusion equation. This is due to 
the fact that in ever smaller devices the assumption that the
free carriers can be modeled as a continuum becomes invalid. On the
other hand, the drift diffusion equations are derived by a scaling limit
process, where the mean free path of a particle tends to zero.

\paragraph{The inverse problems}

This paper is devoted to the investigation of inverse problems related
to {\em drift diffusion equations} modeling semiconductor devices.
In this context we analyze several inverse problems related to the
identification of doping profiles. In all these inverse problems the
parameter to be identified corresponds to the so called {\em doping profile}. 
Such profile enters as a 
functional parameter in a system of PDE's. However, the reconstruction
problems are related to data generated by different types of measurement
techniques.

Identification problems for semiconductor devices, although of increasing 
technological importance, seem to be poorly understood so far.
In the inverse problem literature there has been increasing interest 
on the identification of a position dependent function $C = C(x)$
representing the doping profile, i.e., the density difference of ionized
donors and acceptors. These are the so-called {\em inverse doping profile
problems}. See, for example, \cite{BELM04,BEMP01,BEM02, LMZ05, FI92, FIR02, FI94, BFI93}
and references therein. 

In some cases, e.g., the p-n diode, it may be assumed that the function $C$
is piecewise constant over the device. In this case, the problem reduces to
identifying the curves (or surfaces) between the subdomains where doping is
constant. Particularly important are the curves separating subdomains where
the doping profile assumes constant values of different signs. These curves
are called {\em pn-junctions} (see Section~\ref{sec:2} for details).
In the {\em ion implantation} technique, the most important technique of
silicon devices, only a rough estimate of the doping profile can be obtained
by process modelling (see, e.g., \cite{Sel84}). An efficient alternative to
determine the real doping profile is the use of reconstruction methods
from indirect data.

Another relevant inverse problem concerns identifying transistor contact
resistivity of planar electronic devices, such as MOSFETs (metal oxide
semiconductor field-effect transistors) is treated in \cite{FC92}. It is shown
that a one-point boundary measurement of the potential is sufficient to
identify the resistivity from a one-parameter monotone family, and such
identification is both stable and continuously dependent on the
parameter. Because of the device miniaturization, it is impossible to measure
the contact resistivity in a direct way to satisfactory accuracy. There are
extensive experimental and simulation studies for the determination of contact
resistivity by certain accessible boundary measurements.

Yet another inverse problem is that of determining the contact resistivity
of a semiconductor device from a single voltage measurement \cite{BF91}.
It can be modeled as an inverse problem for the elliptic
differential equation $\Delta V - p \chi(S)u = 0$ in $\Omega \subset \mathbb
R^2$, $\partial V / \partial n = g \geq 0$ but $g \not\equiv 0$ on
$\partial\Omega$, where $V(x)$ is the measured voltage, $S \subset \Omega$ and
$p>0$ are unknown. In \cite{BF91}, the authors consider the identification of
$p$ when the contact location $S$ is also known.

\paragraph{Outline of the article}

In Section~\ref{sec:2} we introduce and discuss relevant properties of
the main mathematical models: the (transient and stationary) systems of
drift diffusion equations. \\
In Section~\ref{sec:3} we derive, from the drift diffusion equations, some
special stationary and transient models, which will serve as mathematical
background to the formulation the inverse doping profile problems. \\
In Section~\ref{sec:4} we formulate several inverse problems, which relates
to specific measurement procedures for the voltage-current map (namely
{\em pointwise measurements of the current density} and {\em current flow
measurements through a contact}) as well as to specific model idealizations. \\
In Section~\ref{sec:5} we present a short description of techniques from 
the theory of inverse problems that is
used to handle the doping profile identification problem described in the other sections. \\
In Section~\ref{sec:6} we present numerical experiments for some models
concerning the inverse doping profile problem for the stationary linearized
unipolar and bipolar cases.

\section{Drift diffusion equations} \label{sec:2}

\paragraph{The transient model}

The basic semiconductor device equations in the {\em transient case}
consist of the Poisson equation (\ref{eq:VanR1}), the continuity equations
for electrons (\ref{eq:VanR2}) and holes (\ref{eq:VanR3}), and the current
relations for electrons (\ref{eq:VanR4}) and holes (\ref{eq:VanR5}). For
some applications, in order to account for thermal effects in semiconductor
devices, its also necessary to add to this system the heat flow equation
(\ref{eq:VanR6}).
\begin{subequations} \label{eq:VanR} \begin{eqnarray}
{\rm div} (\epsilon \nabla V) &       \hskip-1.2cm \label{eq:VanR1}
   =  q(n - p - C) & {\rm in}\ \Omega \times (0,T) \\
{\rm div}\, J_n &                     \hskip-1.5cm \label{eq:VanR2}
   =  q ( \partial_t n + R) & {\rm in}\ \Omega \times (0,T) \\
{\rm div}\, J_p &                     \hskip-1.25cm \label{eq:VanR3}
   =  q (-\partial_t p - R) & {\rm in}\ \Omega \times (0,T)  \\
J_n &                                 \hskip-0.15cm \label{eq:VanR4}
   =  q ( D_n \nabla n - \mu_n n \nabla V) & {\rm in}\ \Omega \times (0,T) \\
J_p &                                 \label{eq:VanR5}
   =  q (-D_p \nabla p - \mu_p p \nabla V) & {\rm in}\ \Omega \times (0,T)  \\
\rho\ c\ \partial_t \mathcal T - H &  \hskip-1.3cm \label{eq:VanR6}
   = {\rm div}\,k(\mathcal T) \nabla \mathcal T
   & {\rm in}\ \Omega \times (0,T) \, .
\end{eqnarray} \end{subequations}
This system is defined in $\Omega \times (0,T)$, where $\Omega \subset
\mathbb{R}^d$ ($d=1,2,3$) is a domain representing the semiconductor
device.
Here $V$ denotes the electrostatic potential ($- \nabla V$ is the electric
field $E$), $n$ and $p$ are the concentration of free
carriers of negative charge (electrons) and positive charge (holes)
respectively and $J_n$ and $J_p$ are the densities of the electron and the
hole current respectively. $D_n$ and $D_p$ are the diffusion coefficients for
electrons and holes respectively. $\mu_n$ and $\mu_p$ denote the mobilities of
electrons and holes respectively. The positive constants $\epsilon$ and $q$
denote the permittivity coefficient (for silicon) and the elementary charge.

The function $R$ has the form $R = \mathcal{R}(n,p,x) (np - n_i^2)$ and
denotes the {\em recombination-generation rate} ($n_i$ is the intrinsic
carrier density). The {\em bandgap} is relatively large for semiconductors
(gap between valence and conduction band), and a significant amount of energy
is necessary to transfer electrons from the valence and to the conduction
band. This process is called generation of electron-hole pairs. On the
other hand, the reverse process corresponds to the transfer of a
conduction electron into the lower energetic valence band. This process
is called recombination of electron-hole pairs. In our model these
phenomena are described by the recombination-generation rate $R$.
Frequently adopted in the literature are the Shockley-Read-Hall model
($\mathcal{R}_{SRH}$) and the Auger model ($\mathcal{R}_{AU}$). They are defined by
$$  \mathcal{R}_{SRH} \bydef \frac{1}{\tau_p(n+n_i) + \tau_p(p+n_i)}
    \, , \ \ \mathcal{R}_{AU} \bydef (C_n n + C_p p)  \, ,$$
where $C_n$, $C_p$, $\tau_n$ and $\tau_p$ are positive constants
whose physical values are listed in the Appendix.

The function $\mathcal T$ represents the temperature and the constants $\rho$
and $c$ denote the specific mass density and specific heat of the material
respectively. Furthermore, $k(T)$ and $H$ denote the thermal conductivity and the locally
generated heat. Equation (\ref{eq:VanR6}) was presented here only for the
sake of completeness of the model and shall not be considered in the
subsequent development.

The function $C(x)$ models a preconcentration of ions in the crystal, so $C(x)
= C_{+}(x) - C_{-}(x)$ holds, where $C_{+}$ and $C_{-}$ are concentrations of
negative and positive ions respectively. In those subregions of $\Omega$ for
which the preconcentration of negative ions predominate (P-regions), we have
$C(x) < 0$. Analogously, we define the N-regions, where $C(x) > 0$ holds.
The boundaries between the P-regions and N-regions (where $C$ change sign)
are called {\em pn-junctions}.

In the sequel we turn our attention to the boundary conditions.
We assume the boundary $\partial\Omega$ of $\Omega$ to be divided into two
nonempty disjoint parts: $\partial\Omega = \overline{\partial\Omega_N} \cup
\overline{\partial\Omega_D}$. The Dirichlet part of the boundary
$\partial\Omega_D$ models the Ohmic contacts, where the potential $V$ as
well as the concentrations $n$ and $p$ are prescribed. The Neumann part
$\partial\Omega_N$ of the boundary corresponds to insulating surfaces,
thus a zero current flow and a zero electric field in the normal
direction are prescribed.
The Neumann boundary conditions for system (\ref{eq:VanR1}) -- (\ref{eq:VanR5})
read:
\begin{equation} \label{eq:VanR-bcN}
 \frac{\partial V}{\partial\nu}(x,t) =
   \frac{\partial n}{\partial\nu}(x,t) =
   \frac{\partial p}{\partial\nu}(x,t) = 0 \, , \
   \partial\Omega_N \times [0,T] \, .
\end{equation}
Moreover, at $\partial\Omega_D \times [0,T]$, the following Dirichlet boundary
conditions are imposed:
\begin{subequations} \label{eq:VanR-bcD} \begin{eqnarray}
V(x,t) & = & V_D(x,t)
       \, = \, U(x,t) + V_{\rm bi}(x)
       \, = \, U(x,t) + U_T \, \ln (n_D(x)/ n_i)  \\[1ex]
n(x,t) & = & n_D(x)
       \, = \, \frac{1}{2} \left(C(x) + \sqrt{C(x)^2 + 4 n_i^2}\right) \\[1ex]
p(x,t) & = & p_D(x)
       \, = \, \frac{1}{2} \left(-C(x) + \sqrt{C(x)^2 + 4 n_i^2}\right) \, .
\end{eqnarray} \end{subequations}
Here, the function $U(x,t)$ denotes the applied potential.
We shall consider the simple situation $\partial\Omega_D = \Gamma_0 \cup
\Gamma_1$, which occurs, e.g., in a diode. The disjoint boundary parts
$\Gamma_i$, $i=0,1$, correspond to distinct contacts.
Differences in $U(x)$ between different segments of $\partial\Omega_D$
correspond to the applied bias between these two contacts.
The constant $U_T$ represents the thermal voltage. Moreover the initial
conditions $n(x,0) \ge 0$, $p(x,0) \ge 0$ have to be imposed.

We conclude this paragraph discussing the solution theory for the transient
drift diffusion system (\ref{eq:VanR1}) -- (\ref{eq:VanR5}),
(\ref{eq:VanR-bcN}), (\ref{eq:VanR-bcD}).
\footnote{In order to simplify th model we neglect thermal effects.}
Twenty years ago, existence and uniqueness of global in time solutions for
the transient drift diffusion equations were demonstrated by Gajewski in
\cite{Gaj85}. Under the assumption that the doping profile satisfies,
$C \in L^r(\Omega)$ for $d \le r \le 6$, it is shown that
\begin{multline} \label{eq:regul-tdde}
(V - V_D, n - n_D, p - p_D) \in W \bydef \\
\big\{ C( [0,T]; H^2_0(\Omega) ) \cap L^2( [0,T]; W^{2,r}_0(\Omega) )
       \cap H^1( [0,T]; \widetilde W ) \big\}
 \times \widetilde{\widetilde W} \times \widetilde{\widetilde W} \, ,
\end{multline}
where $\widetilde W \bydef \{ w \in H^1(\Omega) ; \ w|_{\partial\Omega_D} = 0 \}$
and
$$ \widetilde{\widetilde W} \bydef C( [0,T]; L^2(\Omega) ) \cap
   L^2( [0,T]; \widetilde W) \cap H^1( [0,T]; \widetilde W^*) \, . $$

In the special one dimensional case $\Omega = (0,L)$, a stronger result
is proved, namely:

\begin{lemma}
Let the doping profile satisfy $C \in L^r(\Omega)$, for $d \le r \le 6$.
If the mobilities $\mu_n$ and $\mu_p$ are in $L^\infty(\Omega)$, 
then every solution $(V,n,p)$ of the transient drift diffusion equations
(\ref{eq:VanR1}) -- (\ref{eq:VanR5}), (\ref{eq:VanR-bcN}) and
(\ref{eq:VanR-bcD}) satisfies (\ref{eq:regul-tdde}). Moreover,
$$ (V,n,p) \in C([0,T]; H^2(\Omega)) \cap C([0,T]; W^{1,\infty}(\Omega))^2 . $$
\end{lemma}

\paragraph{The stationary model}

In this paragraph we turn our attention to the stationary drift diffusion
equations. We neglect the thermal effects and assume further
$\frac{\partial n}{\partial t} = \frac{\partial p}{\partial t} = 0$.
Thus, the {\em stationary drift diffusion model} is derived from
(\ref{eq:VanR1}) -- (\ref{eq:VanR5}) in a straightforward way.
Next, motivated by the Einstein relations $D_n = U_T \mu_n$ and
$D_p = U_T \mu_p$ (a standard assumption about the mobilities and
diffusion coefficients), one introduces the so-called
{\em Slotboom variables} $u$ and $v$. They are related to the
original $n$ and $p$ variables by the formula:
\begin{equation} \label{eq:slotboom}
n(x) = n_i \exp\left(\frac{ V(x)}{U_T}\right)\, u(x)\, , \ \ \
p(x) = n_i \exp\left(\frac{-V(x)}{U_T}\right)\, v(x) \, .
\end{equation}
For convenience, we rescale the potential and the mobilities, i.e.
$ V(x) \ \leftarrow \ V(x) / U_T$, \,$\mu_n \leftarrow q U_T \mu_n$,\,
$\mu_p \leftarrow q U_T \mu_p$. It is obvious to check that the current
relations now read $J_n = \mu_n n_i \, e^{ V} \nabla u$,\,
$J_p = -\mu_p n_i \, e^{-V} \nabla v$.

Now we can write the stationary drift diffusion equations in the form
\begin{subequations} \label{eq:dd-sys-nlin} \begin{eqnarray}
\lambda^2 \, \Delta V & \hskip-0.2cm \label{eq:dd-sys1}
  = \ \delta^2 \big(e^Vu - e^{-V}v\big) - C(x) & {\rm in}\ \Omega \\
{\rm div}\, J_n & \hskip-0.3cm \label{eq:dd-sys2}
  = \ \delta^4 \, Q(V,u,v,x) \, (u v - 1)      & {\rm in}\ \Omega \\
{\rm div}\, J_p & \label{eq:dd-sys3}
  = \ - \delta^4 \, Q(V,u,v,x) \, (u v - 1)    & {\rm in}\ \Omega \\[1ex]
V & \hskip-1.45cm = \ V_D \ = \ U + V_{\rm bi} & {\rm on}\ \partial\Omega_D
\label{eq:dd-sys4} \\
u & \hskip-2.0cm = \ u_D \ = \ e^{-U}          & {\rm on}\ \partial\Omega_D
\label{eq:dd-sys5} \\
v & \hskip-2.3cm = \ v_D \ = \ e^{U}           & {\rm on}\ \partial\Omega_D
\label{eq:dd-sys6} \\
\nabla V \cdot \nu & \hskip-0.6cm = \ J_n\cdot\nu \ = \ J_p\cdot\nu \ = \ 0
                   & \rm on \ \partial\Omega_N \, , \label{eq:dd-sys7}
%
\end{eqnarray} \end{subequations}
where $\lambda^2 \bydef \epsilon/(q U_T)$ is the Debye length of the device,
$\delta^2 \bydef n_i$ and the function $Q$ is defined implicitly by the
relation $Q(V,u,v,x) = {\cal R}(n,p,x)$.
\footnote{Notice the applied potential has also to be rescaled:
$U(x) \leftarrow U(x) / U_T$.}

One should notice that, due to the thermal equilibrium assumption, it follows
$np = n_i^2$, and the assumption of vanishing space charge density gives
$n-p-C = 0$, for $x \in \partial\Omega_D$. This fact motivates the
boundary conditions on the Dirichlet part of the boundary.

It is worth mentioning that, in a realistic model, the mobilities
$\mu_n$ and $\mu_p$ usually depend on the electric field strength
$|\nabla V|$.
In what follows, we assume that $\mu_n$ and $\mu_p$ are positive constants.
This assumption simplifies the subsequent analysis, allowing us to concentrate
on the inverse doping problems. As a matter of fact, this dependence could
be incorporated in the model without changing the results described in the
sequel.

Next, we describe some existence and uniqueness results for the 
stationary drift diffusion equations. We start presenting a classical
existence result

\begin{lemma}{\bf \cite[Theorem 3.3.16]{MRS90}} \label{prop:MRS-3316}
Let $\kappa > 1$ be a constant satisfying $\kappa^{-1} \le u_D(x)$,
$v_D(x) \le \kappa$, $x \in \partial\Omega_D$, and let
$-\infty < C_m \leq C_M < + \infty$. Then for any $C \in
\{ L^\infty(\Omega); \ C_m \le C(x) \le C_M, \ x \in \Omega \}$,
the system (\ref{eq:dd-sys1}) -- (\ref{eq:dd-sys7}) admits a weak
solution $(V,u,v) \in ( H^1(\Omega) \cap L^\infty(\Omega) )^3$.
%
%
\end{lemma}

Under stronger assumptions on the boundary parts $\partial\Omega_D$,
$\partial\Omega_N$ as well as on the boundary conditions $V_D$, $u_D$,
$v_D$, it is even possible to show $H^2$-regularity for a solution $(V,u,v)$
of system (\ref{eq:dd-sys1}) -- (\ref{eq:dd-sys7}). For details on this
result we refer the reader to \cite[Theorem 3.3.1]{MRS90}.

As far as uniqueness of solutions of system
(\ref{eq:dd-sys1}) -- (\ref{eq:dd-sys7}) is concerned, some results can
be obtained if the applied voltage is small (in the norm of
$L^\infty(\partial\Omega_D) \cap H^{3/2}(\partial\Omega_D)$).

\begin{lemma}{\bf \cite[Theorem 2.4]{BEMP01}} \label{prop:BEMP-24}
Let the applied voltage $U$ be such that $\|U\|_{L^\infty(\partial\Omega_D)}
+ \|U\|_{H^{3/2}(\partial\Omega_D)}$ is sufficiently small. Then, system
(\ref{eq:dd-sys1}) -- (\ref{eq:dd-sys7}) has a unique solution $(V,u,v) \in
( H^1(\Omega) \cap L^\infty(\Omega) )^3$.
\end{lemma}

Since existence and uniqueness of solutions for system
(\ref{eq:dd-sys1}) -- (\ref{eq:dd-sys7}) can only be guaranteed for
small applied voltages, it is reasonable to consider, instead of
this system, its linearized version around the equilibrium point
$U \equiv 0$. We shall return to this point in the next section, where
the voltage-current map is introduced.

\section{Special models} \label{sec:3}

In the next subsections we assume several different simplifications of
the drift diffusion models introduced in Section~\ref{sec:2} and
derive some special cases which will serve as underlying models for the
inverse problems investigated in Section~\ref{sec:4}.

\subsection{The linearized stationary drift diffusion equations
(close to equilibrium)} \label{ssec:stat-lin-equil}

We begin this subsection by introducing the {\em thermal equilibrium}
assumption for the stationary drift diffusion equations. This is a previous
step to derive a linearized system of stationary drift diffusion equations
(close to equilibrium).

The thermal equilibrium assumption refers to the condition in which the
semiconductor is not subject to external excitations, except for a uniform
temperature, i.e. no voltages or electric fields are applied.
It is worth noticing that, under the thermal equilibrium assumption, all
externally applied potentials to the semiconductor contacts are zero
(i.e. $U(x) = 0$). Moreover, the thermal generation is perfectly balanced
by recombination (i.e. $\mathcal R = 0$).

If the applied voltage satisfies $U = 0$, one immediately sees that the
solution of system (\ref{eq:dd-sys1}) -- (\ref{eq:dd-sys7}) simplifies to
$(V,u,v) = (V^0, 1,1)$, where $V^0$ solves
\begin{subequations}  \label{eq:equil-case} \begin{eqnarray}
\lambda^2 \, \Delta V^0 & \label{eq:equil-caseA}
 = \ e^{V^0} - e^{-V^0} - C(x)                & {\rm in}\ \Omega \\
V^0 & \hskip-2.0cm
 = \ V_{\rm bi}(x)                            & {\rm on}\ \partial\Omega_D \\
   \nabla V^0 \cdot \nu & \hskip-2.7cm
 = \ 0                                        & {\rm on}\ \partial\Omega_N \, .
\end{eqnarray} \end{subequations}

For some of the models discussed below, we will be interested in the linearized
drift diffusion system at the equilibrium. Keeping this in mind, we compute
the Gateaux derivative of the solution of system (\ref{eq:dd-sys1}) --
(\ref{eq:dd-sys7}) with respect to the voltage $U$ at the point $U \equiv 0$
in the direction $h$. This directional derivative is given by the solution
$(\hat V, \hat u, \hat v)$ of
\begin{subequations}  \label{eq:dd-sys-lin} \begin{eqnarray}
 \lambda^2 \, \Delta \hat V & 
 = \ e^{V^0} \hat u + e^{-V^0} \hat v + ( e^{V^0} + e^{-V^0} ) \hat V
                                                & {\rm in}\ \Omega \\
   {\rm div}\, (\mu_n e^{V^0} \nabla \hat u) & \hskip-2.25cm
 = \ Q_0(V^0,x) (\hat u + \hat v)               & {\rm in}\ \Omega \\
   {\rm div}\, (\mu_p e^{-V^0} \nabla \hat v) & \hskip-2.21cm
 = \ Q_0(V^0,x) (\hat u + \hat v)               & {\rm in}\ \Omega \\[1ex]
   \hat V & \hskip-4.6cm
 = \ h                                          & {\rm on}\ \partial\Omega_D \\
   \hat u & \hskip-4.3cm
 = \ -h                                         & {\rm on}\ \partial\Omega_D \\
   \hat v & \hskip-4.6cm
 = \  h                                         & {\rm on}\ \partial\Omega_D \\
   \nabla V^0 \cdot\nu & \hskip-1.5cm
 = \ \nabla\hat u \cdot\nu \, = \, \nabla\hat v \cdot\nu \, = \, 0
                                              & {\rm on}\ \partial\Omega_N \, ,
\end{eqnarray} \end{subequations}
where the function $Q_0$ satisfies $Q_0(V^0,x) = Q(V^0,1,1,x)$.

\subsection{Linearized stationary bipolar case (close to equilibrium)}
\label{ssec:bipol}

In this subsection we present a special case, which plays a key r\^ole in 
modeling inverse doping problems related to {\em current-flow} measurements.

The discussion is motivated by the {\em stationary voltage-current map}
(V-C) map
$$ \begin{array}{rcl}
   \Sigma_C: H^{3/2}(\partial\Omega_D) & \to & \mathbb R \, . \\
   U & \mapsto & \displaystyle\int_{\Gamma_1} (J_n+J_p)\cdot\nu \, ds
   \end{array} $$
Here $(V,u,v)$ is the solution of (\ref{eq:dd-sys-nlin}) for an applied
voltage $U$. This operator models practical experiments where
{\em voltage-current data} are available, i.e. measurements of the
averaged outflow current density on $\Gamma_1 \subset \partial\Omega_D$.

The {\em linearized stationary bipolar case (close to equilibrium)}
corresponds to the model obtained from the drift diffusion equations
(\ref{eq:dd-sys-nlin}) by linearizing the V-C map at $U \equiv 0$. This
simplification is motivated by the fact that, due to hysteresis effects for
large applied voltage, the V-C map can only be defined as a single-valued
function in a neighborhood of $U=0$.
Moreover, the following simplifying assumptions are also taken into account:

\begin{itemize}
\item[{\it A1)}] \ The electron mobility $\mu_n$ and hole mobility $\mu_p$
are constant;
\item[{\it A2)}] \ No recombination-generation rate is present, i.e.
${\cal R} = 0$ (or $Q_0 = 0$).
\end{itemize}

An immediate consequence of our assumptions is the fact that the 
Poisson equation and the continuity equations decouple. Indeed,
from (\ref{eq:dd-sys-lin}) we see that the Gateaux derivative of
the V-C map $\Sigma_C$ at the point $U=0$ in the direction
$h \in H^{3/2}(\partial\Omega_D)$ is given by the expression
\begin{equation} \label{eq:def-sigma-prime-C}
\Sigma'_C(0) h = \int_{\Gamma_1}
\left( \mu_n \, e^{V_{\rm bi}} \hat{u}_\nu -
\mu_p \, e^{-V_{\rm bi}} \hat{v}_\nu \right) \, ds ,
\end{equation}
where $(\hat{u}, \hat{v})$ solve
\begin{subequations}  \label{eq:bipol-stat}  \begin{eqnarray}
{\rm div}\, (\mu_n e^{V^0} \nabla \hat{u})   & \hskip-1.6cm
 = \ 0                \label{eq:bipol-statA} & {\rm in}\ \Omega \\
{\rm div}\, (\mu_p e^{-V^0} \nabla \hat{v})  & \hskip-1.6cm
 = \ 0                \label{eq:bipol-statB} & {\rm in}\ \Omega \\[1ex]
\hat{u} & \hskip-1.3cm
 = \ -h                                      & {\rm on}\ \partial\Omega_D \\
\hat{v} & \hskip-1.55cm
 = \ h                                       & {\rm on}\ \partial\Omega_D \\
\nabla\hat{u} \cdot\nu &
 = \ \nabla\hat{v} \cdot\nu \ = \ 0          & {\rm on}\ \partial\Omega_N
\end{eqnarray} \end{subequations}
and $V^0$ is the solution of the equilibrium problem (\ref{eq:equil-case});
(see Lemma~\ref{prop:bemp31} for details).

Notice that the solution of the Poisson equation can be computed a priori,
since it does not depend on $h$. The linear operator $\Sigma'_C(0)$
is continuous. Actually, we can prove more: since $(u,v)$ depend
continuously in $H^2(\Omega)^2$ on the boundary data $h$ in
$H^{3/2} (\partial\Omega_D)$, it follows from the boundedness and
compactness of the trace operator $\gamma: H^2(\Omega) \to
H^{1/2}(\Gamma_1)$ that $\Sigma'_C(0)$ is a compact operator. The
operator $\Sigma_C'(0)$ maps the Dirichlet data for $(\hat{u}, \hat{v})$
to a weighted sum of their Neumann data and can be compared with the
DtN map in the electrical impedance tomography. See \cite{Bor02,BU02,Nac96}.

\subsection{Linearized stationary unipolar case (close to equilibrium)}
\label{ssec:upol}

The linearized unipolar case (close to equilibrium) corresponds to the model
obtained from the unipolar drift diffusion equations by linearizing the V-C
map at $U \equiv 0$. Additionally to {\it A1)} and {\it A2)}, we further assume:

\begin{itemize}
\item[{\it A3)}] \ The concentration of holes satisfy $p = 0$ (or,
equivalently, $v = 0$ in $\Omega$).
\end{itemize}

Under those assumptions, the Gateaux derivative of the V-C map $\Sigma_C$
at the point $U=0$ in the direction $h$ is given by
$$ \Sigma'_C(0) h = \int_{\Gamma_1}
   \mu_n \, e^{V_{\rm bi}} \hat{u}_\nu \, ds $$
where $\hat{u}$ solve
\begin{subequations} \label{eq:unipol-stat} \begin{eqnarray}
{\rm div}\, (\mu_n e^{V^0} \nabla \hat{u}) & \hskip-0.6cm
 = \ 0                                             & {\rm in}\ \Omega \\
\hat{u} & = \ - h(x)                               & {\rm on}\ \Omega_D \\
\nabla\hat{u} \cdot\nu & \hskip-0.6cm = \ 0        & {\rm on}\ \Omega_N
\end{eqnarray} \end{subequations}
and $V^0$ is the solution of the equilibrium problem (\ref{eq:equil-case}),
with (\ref{eq:equil-caseA}) replaced by
\begin{center}
\hfill
\mbox{$\lambda^2 \, \Delta V^0 = \ e^{V^0} - C(x)$ \ in\ $\Omega$.}
\hfill (\ref{eq:equil-caseA}')
\end{center}

\subsection{Linearized transient bipolar case (close to equilibrium)}
\label{ssec:bipol-trans}

In this subsection we introduce a transient case, which is the time dependent
counterpart of the bipolar model discussed in Subsection~\ref{ssec:bipol}. It
will serve as background for the formulation of inverse doping problems
related to transient current flow measurements.

As in Subsection~\ref{ssec:bipol}, we begin the discussion by introducing
the {\em transient voltage-current map}. For an applied time dependent
voltage $U(x,t)$, the transient V-C map is given by
\begin{equation} \label{eq:def-sigma-tC}
\begin{array}{rcl}
   \Sigma_{t,C}: L^2([0,T]; H^{3/2}(\partial\Omega_D)) & \to & L^2(0,T) \, . \\
   U(\cdot,t) & \mapsto & \displaystyle\int_{\Gamma_1}
                          [ J_n(\cdot,t) + J_p(\cdot,t) ] \cdot\nu \, ds
   \end{array}
\end{equation}
Here $(V,n,p)$ is the solution of (\ref{eq:VanR}), (\ref{eq:VanR-bcN}),
(\ref{eq:VanR-bcD}) for an applied voltage $U$.%
\footnote{Once more we disconsider equation (\ref{eq:VanR6}).}
This operator models practical experiments where time dependent
{\em voltage-current data} are available.
In \cite{BEM02} it is shown that the nonlinear operator $\Sigma_{t,C}$ is
well defined, continuous and Fr\'echet differentiable. In the sequel we
derive the Gateaux derivative of $\Sigma_{t,C}$ in equilibrium.

As in the stationary cases, we shall consider the transient drift diffusion
equations under the {\em thermal equilibrium} assumption. It is immediate
to observe that, for zero applied voltage $U(\cdot,t) = 0$, the solution
$(V^0,n^0,p^0)$ of (\ref{eq:VanR}), (\ref{eq:VanR-bcN}), (\ref{eq:VanR-bcD})
is constant in time, being the counterpart (in the $n$, $p$ variables)
of the solution triplet $(V^0,1,1)$ in the slotboom variables
(see (\ref{eq:slotboom})).

Here again, we assume $A1)$, $A2)$ of Subsection~\ref{ssec:bipol}.
Then, arguing as in Subsection~\ref{ssec:stat-lin-equil}, it follows that the
Gateaux derivative of the transient V-C map $\Sigma_{t,C}$ at the point
$U = 0$ in the direction $h(\cdot,t) \in L^2( [0,T] ; 
H^{3/2}(\partial\Omega_D))$ is given by
\begin{equation} \label{eq:def-sigma-prime-tC}
\Sigma'_{t,C}(0) h = \int_{\Gamma_1}
\big[ \mu_n ( {\hat n}_\nu - {\hat n} V^0_\nu - n^0 {\hat V}_\nu )
     -\mu_p ( {\hat p}_\nu + {\hat p} V^0_\nu + p^0 {\hat V}_\nu ) \big]\, ds ,
\end{equation}
where $(\hat V, \hat n, \hat p)$ solve
\begin{subequations} \label{eq:VanR-lin-tr} \begin{eqnarray}
\lambda^2 \hat V &          \hskip-3.80cm      \label{eq:VanR-lin-trA}
   =  \hat n - \hat p       & {\rm in}\ \Omega \times (0,T) \\
\partial_t \hat n &                            \label{eq:VanR-lin-trB}
   =  {\rm div} ( \mu_n [\nabla\hat n - \hat n \nabla V^0 - n^0 \nabla\hat V] )
                            & {\rm in}\ \Omega \times (0,T) \\
\partial_t \hat p &         \hskip-0.13cm      \label{eq:VanR-lin-trC}
   =  {\rm div} ( \mu_p [\nabla\hat p + \hat p \nabla V^0 + p^0 \nabla\hat V] )
                            & {\rm in}\ \Omega \times (0,T)  \\
\hat V &                    \hskip-4.40cm      \label{eq:VanR-lin-trD}
   = h                      & {\rm on}\ \partial\Omega_D \times (0,T) \\
\hat n &                    \hskip-3.80cm      \label{eq:VanR-lin-trE}
   = \hat p = 0             & {\rm on}\ \partial\Omega_D \times (0,T) \\
\nabla\hat{V} \cdot\nu &    \hskip-1.62cm      \label{eq:VanR-lin-trF}
   = \nabla\hat{n} \cdot\nu = \nabla\hat{p} \cdot\nu = 0 
                            & {\rm on}\ \partial\Omega_N \times (0,T) \, .
\end{eqnarray} \end{subequations}

Notice that, differently from the stationary case, the Poisson equation
(\ref{eq:VanR-lin-trA}) and the continuity equations (\ref{eq:VanR-lin-trB}),
(\ref{eq:VanR-lin-trC}) do not decouple.

\section{Inverse Problems for Semiconductors} \label{sec:4}

In practical experiments there are different types of measurement
techniques, such as
\begin{itemize}
\item {\em Laser-beam-induced-current} (LBIC) measurements;
\item {\em Capacitance} measurements;
\item {\em Current Flow} measurements.
\end{itemize}

We refer to \cite{FI92,FI94,FIR02} for the first type and to
\cite{BELM04,BEMP01,BEM02,LMZ05} for the last two types.
These measurement techniques are related to different types of
data and lead to different inverse problems for reconstructing
the doping profile. They are the so-called {\em inverse doping
profile problems}.
In the following subsections we address inverse problems related
to each one of these measurement techniques.

\subsection{The stationary voltage-current map} \label{ssec:ip-vcmap}

We begin this subsection verifying that the voltage-current map
$\Sigma_C$, introduced in Subsection~\ref{ssec:bipol}, is well defined
in a suitable neighborhood of $U=0$.

\begin{lemma} \mbox{\bf \cite[Proposition 3.1]{BEMP01}}
\label{prop:bemp31}
For each applied voltage $U \in B_r(0) \subset H^{3/2}(\partial\Omega_D)$
with $r>0$ sufficiently small, the current $J \cdot \nu \in H^{1/2}
(\Gamma_1)$ is uniquely defined. Furthermore, $\Sigma_C:
H^{3/2}(\partial\Omega_D) \to H^{1/2}(\Gamma_1)$ is continuous
and is continuously differentiable in $B_r(0)$. Moreover, its derivative
in the direction $h \in H^{3/2}(\partial\Omega_D)$ is given by the operator
$\Sigma'_C(0)$ defined in (\ref{eq:def-sigma-prime-C}).
\end{lemma}

As a matter of fact, we can actually prove that, since $(\hat u,\hat v)$
in (\ref{eq:bipol-stat}) depend continuously (in $H^2(\Omega)^2$) on the
boundary data $U \in H^{3/2}(\partial\Omega_D)$, it follows from the
boundedness and compactness of the trace operator $\gamma: H^2(\Omega)
\to H^{1/2}(\Gamma_1)$ that $\Sigma'_C(0)$ is a bounded and compact
operator. The operator $\Sigma_C'(0)$ in  (\ref{eq:def-sigma-prime-C})
maps the Dirichlet data for $(\hat u,\hat v)$ to a weighted sum of their
Neumann data and the related inverse problem can be compared with the
identification problem in {\em Electrical Impedance Tomography} (EIT).

Proposition~\ref{prop:bemp31} establishes a basic property to consider
the inverse problem of reconstructing the doping profile $C$ from the
V-C map. In the sequel we shall consider two possible inverse problems
for the V-C map.

\paragraph{Current flow measurements through a contact}

In the first inverse problem we assume that, for each $C$, the output
is  given by $\Sigma_C'(0) U_j$ for some $U_j$. A realistic experiment
corresponds to measure, for given $\{ U_j \}_{j=1}^N$, with $\|U_j\|$
small, the outputs
$$ \big\{ \Sigma_C'(0) U_j\ |\ \ j=1,\cdots,N \big\} $$
(recall that $\Sigma_C(0) = (V^0,1,1)$).
In practice, the functions $U_j$ are chosen to be piecewise constant on
the contact $\Gamma_1$ and to vanish on $\Gamma_0$. 
From the definition of $\Sigma_C'(0)$ we deduce the following abstract
formulation of the inverse doping profile problem for the V-C map:
\begin{equation} \label{eq:ip-abstract}
 F(C) \ = \ Y \, ,
\end{equation}
where
\begin{enumerate}
\item[1)] $\{ U_j \}_{j=1}^N \subset H^{3/2}(\partial\Omega_D)$ are
fixed voltage profiles satisfying $U_j |_{\Gamma_1} = 0$;
\item[2)] Parameter: \ $C = C(x) \ \in \ L^2(\Omega) =: \mathcal X$;
\item[3)] Output: \ $Y = \big\{ \Sigma_C'(0) U_j \big\}_{j=1}^N \in
\mathbb R^N =: \mathcal Y$;
\item[4)] Parameter-to-output map: \ $F: \mathcal X \to \mathcal Y$.
\end{enumerate}
The domain of definition of the operator $F$ is
$$ D(F) \bydef \{ C \in L^\infty(\Omega) ; \, C_m \le C(x) \le C_M,
   \mbox{ a.e. in } \Omega \} \, , $$
where $C_m$ and $C_M$ are suitable positive constants.

This approach is motivated by the fact that, in practical applications,
the V-C map can only be defined in a neighborhood of $U=0$ (due to
hysteresis effects for large  applied voltages). The inverse problem
described above corresponds to the problem of identifying the doping
profile $C$ from the linearized V-C map at $U=0$. See the unipolar and
bipolar cases in Subsections~\ref{ssec:bipol} and~\ref{ssec:upol}.

The nonlinear parameter-to-output operator $F$ is well defined and
Fr\'echet differentiable in its domain of definition $D(F)$. This
assertion follows from standard regularity results in PDE theory
\cite[Propositions~2.2 and~2.3]{BELM04}.

It is worth noticing that the solution of the Poisson equation can be
computed {\em a priori}. The remaining problem (coupled system
(\ref{eq:bipol-stat}) for $(\hat u,\hat v)$) is quite similar to the
problem of EIT (see \cite{Bor02, Isa98}). In this inverse problem the
aim is to identify the conductivity $q = q(x)$ in the equation
$$ -{\rm div}\,(q \nabla u) \ = \ f \ \  {\rm in}\ \Omega \, , $$
from measurements of the {\em Dirichlet-to-Neumann map}, which maps
the applied voltage $u|_{\partial\Omega}$ to the electrical flux
$q u_\nu|_{\partial\Omega}$. The map $\Sigma_C'(0)$ sends the
Dirichlet data for $\hat{u}$ and $\hat{v}$ to the weighted sum of their
Neumann data. It can be seen as the counterpart of electrical impedance
tomography for common conducting materials.

\paragraph{Pointwise measurements of the current density}

In the sequel, we investigate a different formulation of the same inverse
problem related to the V-C map considered above.
Differently from the previous paragraph, we shall assume that the
V-C operator maps the Dirichlet data for $\hat{u}$ and $\hat{v}$ in
(\ref{eq:bipol-stat}) to the sum of their Neumann data, i.e.
$$ \begin{array}{rcl}
     \Sigma_C: H^{3/2}(\partial\Omega_D) & \to & H^{1/2}(\Gamma_1) \\
     U & \mapsto & (J_n + J_p)\cdot\nu |_{\Gamma_1}
   \end{array} $$
where functions $V$, $u$, $v$, $J_n$, $J_p$ and $U$ have the same meaning
as in Subsection~\ref{ssec:bipol}. It is immediate to observe that the
Gateaux derivative of the V-C map $\Sigma_C$ at the point $U=0$ in the
direction $h \in H^{3/2}(\partial\Omega_D)$ is given by
\begin{equation}
\Sigma'_C(0) h = \left( \mu_n \, e^{V_{\rm bi}} \hat{u}_\nu -
\mu_p \, e^{-V_{\rm bi}} \hat{v}_\nu \right) |_{\Gamma_1} \, ,
\end{equation}
where $(\hat{u}, \hat{v})$ solve system (\ref{eq:bipol-stat}). Notice that,
for each voltage profile $U$, the V-C map associates a scalar valued function
defined on $\Gamma_1$. In this case, the outputs $\Sigma'_C(0) U_j$ give much
more information about the parameter $C$ than in the case of current flow
measurements.

Again we can derive an abstract formulation of type (\ref{eq:ip-abstract})
for the inverse doping profile problem for the V-C map with pointwise
measurements of the current density. The only difference to the framework
described in the previous paragraph concerns the definition of the Hilbert
space $Y$, which is now defined by:
\begin{enumerate}
\item[3')] Output: \ $Y = \big\{ \Sigma_C'(0) U_j \big\}_{j=1}^N \in
L^2(\Gamma_1)^N =: \mathcal Y$.
\end{enumerate}
The domain of definition of the operator $F$, remains unaltered.

In Section~\ref{sec:6} we shall consider three numerical implementations
concerning inverse doping problems for the V-C map described above, namely:
\begin{enumerate}
\item[i)]
The stationary linearized unipolar model (close to equilibrium)
with current flow measurements through a contact;
\item[ii)]
The stationary linearized unipolar model (close to equilibrium)
with pointwise measurements of the current density.
\item[iii)]
The stationary linearized bipolar model (close to equilibrium)
with pointwise measurements of the current density.
\end{enumerate}

\subsection{The transient voltage-current map} \label{ssec:ip-vcmap-tr}

In the sequel we shall consider inverse problems for the map $\Sigma_{t,C}$
in (\ref{eq:def-sigma-tC}).
As already observed in Subsection~\ref{ssec:bipol-trans}, this V-C map
is well defined, continuous and Fr\'echet differentiable, its Gateaux
derivative in equilibrium $\Sigma'_{t,C}(0)$ being defined by
(\ref{eq:def-sigma-prime-tC}).

As in Subsection~\ref{ssec:ip-vcmap} we investigate two possible inverse
doping problems for the {\em linearized transient bipolar case (close to
equilibrium)}.

\paragraph{Transient current flow measurements through a contact}

Here we assume that, for each $C$, the output corresponds to
$\Sigma_{t,C}'(0) U_j$ for some prescribed $U_j(\cdot,t)$. The experiment
corresponds to measure, for given $\{ U_j(\cdot,t) \}_{j=1}^N$, with
$\|U_j\|$ small in $L^2([0,T]; H^{3/2}(\partial\Omega_D))$, the (averaged)
currents
$$ \big\{ \Sigma_{t,C}'(0) U_j(\cdot,t)\ |\ \ j=1,\cdots,N \big\} \, . $$
The profile of the voltages $U_j$ is chosen analogously to that of
Subsection~\ref{ssec:ip-vcmap}. Notice that in a realistic transient
experiment, the amplitude of functions $U_j(\cdot,t)$ may vary with the
time, e.g.,
$$  U_j(x,t) \ \bydef \ \left\{ \begin{array}{rl}
      1+t, & |x - x_j| \le h \\
      0  , & {\rm elsewhere} \end{array} \right. $$
where $\Gamma_0 = (0,1) \times \{0\} \subset \mathbb R^2$ and
$0 < x_1 < x_2 < \cdots < x_N < 1$ and $h$ is small enough (compare
with Subsection~\ref{ssec:num-upm}).

The inverse doping profile problem for the V-C map $\Sigma_{t,C}'(0)$ can be
formulated in the abstract form (\ref{eq:ip-abstract}), where
\begin{enumerate}
\item[$1_t$)] $\{ U_j(\cdot,t) \}_{j=1}^N \subset
L^2([0,T]; H^{3/2}(\partial\Omega_D))$ are fixed voltage profiles
satisfying $U_j(\cdot,t) |_{\Gamma_1} = 0$, $t \ge 0$;
\item[$2_t$)] Parameter: \ $C = C(x) \ \in \ L^2(\Omega) =: \mathcal X$;
\item[$3_t$)] Output: \ $Y = \big\{ \Sigma_{t,C}'(0) U_j \big\}_{j=1}^N \in
L^2(0,T)^N =: \mathcal Y$;
\item[$4_t$)] Parameter-to-output map: \ $F: \mathcal X \to \mathcal Y$.
\end{enumerate}
The domain of definition of the operator $F$ is
$$ D(F) \bydef \{ C \in L^\infty(\Omega) ; \, C_m \le C(x) \le C_M,
   \mbox{ a.e. in } \Omega \} \, , $$
where $C_m$ and $C_M$ are suitable positive constants.

From our knowledge about the operator $\Sigma_{t,C}'(0)$ we conclude that
the parameter to output operator $F$ is well defined and continuous. Moreover,
for one dimensional domains $\Omega$ it is shown in \cite{BEM02} that $F$ is
weakly sequentially closed.

\paragraph{Transient pointwise measurements of the current density}

In the previous paragraph, we considered $\Sigma_{t,C}$ to be defined by
(\ref{eq:def-sigma-tC}). Now, we shall assume that, for every time instant
$t \ge 0$, current measurements are available at every point of the segment
$\Gamma_1$. This assumption corresponds to the following definition of the
V-C map:
$$
\begin{array}{rcl}
  \Sigma_{t,C}: L^2([0,T]; H^{3/2}(\partial\Omega_D)) & \to &
                                          L^2([0,T]; H^{1/2}(\Gamma_1)) \, . \\
  U(\cdot,t) & \mapsto & (J_n(\cdot,t) + J_p(\cdot,t)) \cdot\nu
  |_{\Gamma_1} \, ,
\end{array} $$
where $(V,n,p)$ is the solution of (\ref{eq:VanR}), (\ref{eq:VanR-bcN}),
(\ref{eq:VanR-bcD}) for an applied voltage $U(\cdot,t)$. It is immediate
to observe that the Gateaux derivative of $\Sigma_{t,C}$ at the point
$U(\cdot,t) = 0$ in the direction $h \in
L^2([0,T]; H^{3/2}(\partial\Omega_D))$ is given by
$$ \Sigma'_{t,C}(0) h =
[ \mu_n ( {\hat n}_\nu - {\hat n} V^0_\nu - n^0 {\hat V}_\nu )
 -\mu_p ( {\hat p}_\nu + {\hat p} V^0_\nu + p^0 {\hat V}_\nu ) ]
 \big|_{\Gamma_1} \, , $$
where $(\hat V, \hat n, \hat p)$ solve (\ref{eq:VanR-lin-tr}).

The inverse doping profile problem for this V-C map can again be written in
the abstract form $F(C)=Y$. The corresponding framework is now described by
$1_t)$, $2_t)$, $4_t)$ and
\begin{enumerate}
\item[$3'_t$)] Output: \ $Y = \big\{ \Sigma_{t,C}'(0) U_j \big\}_{j=1}^N \in
L^2([0,T]; H^{1/2}(\Gamma_1))^N =: \mathcal Y$;
\end{enumerate}

As in the inverse problem of the previous paragraph, the parameter to output
map $F$ is well defined and continuous. If the domain $\Omega$ is one
dimensional, the results in \cite{BEM02} can be adapted in a straightforward
way and we can conclude that $F$ is weakly sequentially closed.

\section{Background on Inverse Problems and Level Set Methods} \label{sec:5}


In what follows we present some of the tools from the theory of Inverse Problems that are needed
as background to understand the approach we use in the present work. These tools include some classical material
such as for example the singular value decomposition, regularization and Landweber's method, which are treated
in Sections~\ref{svd1}, \ref{reg1}, and \ref{landwebb1}, as well as more recent developments such as use  of
level set methods for handling inverse problems. The latter is treated in Section~\ref{levset}.

\subsection{The Singular Value Decomposition} \label{svd1}
                                                                                
We briefly review the {\em singular value decomposition (SVD)}. This result
has a number of important applications in numerical analysis, inverse
problems, and numerical linear algebra. 

Let $\mathcal{L}(\H,\K)$ denote the space of bounded linear operators from
$\H$ to $\K$, where $\H$ and $\K$ are Hilbert spaces. 
We endow $\mathcal{L}(\H,\K)$ with the {\em uniform operator topology} defined by the norm
$$ ||T||_{\H,\K}
 \bydef \sup_{f\ne 0}  \frac{||T f||_{\K}}{||f||_{\H}} \mbox{ .} $$ 
Whenever no confusion may arise we shall drop the $\H,\K$ subscript in $ ||T||_{\H,\K} $. 
                                                                                                                                              
We recall the concept of compact operator.                                                                                
\begin{definition}
Let $\H$ and $\K$ be Banach spaces  and $T:\H \rightarrow \K$ a linear operator.
$T$ is called {\em compact} (or {\em completely continuous}) if
it maps the unit ball $B_{\H}$ on a pre-compact set, i.e.,
$T(B_{\H})$ has compact closure.
The set of compact operators from $\H$ to $\K$ is denoted by
${\kappa}(\H, \K)$.
\end{definition}

It can be easily shown that the above definition implies that
$T$ is a bounded operator. Furthermore, the set of compact operators
is closed under limits in the uniform operator topology of $\mathcal{L}(\H,\K)$. Typical examples of compact operators are operators of
finite dimensional range.
Another important class of compact operators is given
by integral operators.
 We use the convention that our {\em complex} Hilbert space inner-products
are {\em linear} w. r. t. the first entry and
{\em anti-linear} w. r. t. to the second one.
The following result is instrumental to understand the structure of compact operators.
\begin{theorem} \label{svd}
Let $A \in {\kappa}(\mathcal{H},\mathcal{K}) $, then we can write
\begin{equation}  \label{eqsvd}
A = \sum_{n=1}^r \sigma_n \, (  \cdot \mid  \psi_n  ) \, \phi_n ,
\end{equation}
where $r \in {\mathbb N} \cup \{ \infty \}$ and $\sigma_1 \geq \sigma_2
\geq \cdots \geq \sigma_r > 0 \mbox{ ;}$ the sets $\{ \phi_n \}_{n=1}^r$
and $\{ \psi_n \}_{n=1}^r$ are orthonormal sets (not necessarily
complete) in $\mathcal{K}$ and $ \mathcal{H}$, respectively.
\end{theorem}

As a direct consequence of the SVD decomposition we get that the equation
$A f = g$ has a solution $f$ if, and only if,
\begin{enumerate}
\item The vector $g\in \ker{A^{\ast}}^{\perp}$, and
\item the sum
$$\sum_{n=1}^{r}\frac{1}{\sigma_n^{2}} |(g \mid \phi_{n})|^{2}  < \infty \mbox{ .} $$
\end{enumerate}
In this case the solution of $x$ will be given by
$$ f = \sum_{n=1}^{r}  \frac{1}{\sigma_n} (g \mid \phi_{n}) \psi_n
\mbox{ .} $$
                                                                                
In the finite dimensional case, the transformation $A$ 
will be invertible if, and only if, the dimensions of
${\mathcal{H}}$ and ${\mathcal{K}}$ equal $r$. In other words, the eigenvalues of $A^{\ast} A$  and $A A^{\ast}$ are all
nonzero.
The presence of singular values close
to zero indicates that the solution of the problem $A f = g$
will be doomed to numerical instability.
The condition number of a matrix $A\in {\mathbb{C}}^{n\times n}$ is the
ratio $\sigma_1/\sigma_r$, if $r=n$, and $\infty$ otherwise. See \cite{Golub} for more information
on numerical aspects related to conditioning. 
                                                                                
In the infinite dimensional case, if $r=\infty$ then
the sequence $\left\{ \sigma_n \right\}$ must necessarily converge to zero.
This shows the inherent instability of solving equations of
the form $A f = g$ when $A$ is a compact operator in an infinite
dimensional Hilbert space.

\subsection{Regularization.} \label{reg1}
A mathematical problem, defined in the form of an equation $F(u) = g$ where $F$ is an operator between two Hilbert spaces
$\H$ and $\K$ is said to be well-posed (in the sense of Hadamard) if for every $g\in \K$ the solution $u\in \H$ exists, is unique, and depends continuously on $g$.

Let $A:  \mathcal{H} \rightarrow \mathcal{K}$ be a  compact linear operator.
We now analyze the question of solving a linear equation of the form 
$$ A f = g \mbox{ .} $$ There are three things that can go wrong:
\begin{enumerate}
\item The equation may not be solvable.
(i.e. $g \not \in \mathbf{Ran}(A)$)
\item The solution may not be unique.
(i.e. $A$ is not $1-1$.)
\item The solution may not depend continuously on the data.
(i.e. $A^{-1}$ is not continuous.)
\end{enumerate}
The concept of pseudo-inverse (or generalized inverse) $A^{\dagger}$ is used to handle the cases 1 and 2 above.
We define
$$ A^{\dagger} g = \sum_{n=1}^{r} \frac{1}{\sigma_{n}} ( g | \phi_n ) \psi_n
\mbox{ , } g\in {D}(A^{\dagger}) \mbox{ , } $$ 
where $${D}(A^{\dagger}) \bydef 
\left\{ g \in {\mathcal{K}}  \bigg|  \sum_{n=1}^{r}\frac{1}{\sigma_n^{2}} |(g \mid \phi_{n})|^{2}  < \infty \right\} \mbox{ .} $$
Note that $ A^{\dagger} g$ is the unique solution of $A f = g$ in $(\ker{A})^{\perp}$. 
To tackle the problem of discontinuous $A^{-1}$, one needs to introduce the notion of {\em regularization}.

Let us consider a family of continuous operators 
$T_{\alpha}:\mathcal{K} \rightarrow \mathcal{H}$ such that
\begin{equation}\label{eq1}
\lim_{\alpha\downarrow 0} T_{\alpha} g = A^{\dagger} g \mbox{  , } g \in {D}(A^{\dagger}) \mbox{ .}
\end{equation}
Note that if $A^{\dagger}$ is not bounded, then 
$||T_{\alpha} || \rightarrow \infty $ when $\alpha\downarrow 0$.
If we can solve 
$A f = g $ approximately in the sense that: Let 
$g^{\epsilon}\in \mathcal{K}$ be an approximation to  $g$  such that 
$||g-g^{\epsilon} || \le \epsilon $.
Consider $\alpha(\epsilon)$ such that $\alpha(\epsilon)\downarrow 0$
and $||T_{\alpha(\epsilon)}||\epsilon \rightarrow 0$.
Thus,
\begin{eqnarray*}
|| T_{\alpha(\epsilon)} g^{\epsilon} - A^{\dagger} g ||
 & \le &  ||  T_{\alpha(\epsilon)}(g^{\epsilon} -g) ||   +
||  T_{\alpha(\epsilon)}g - A^{\dagger} g ||   \\
 & \le & ||T_{\alpha(\epsilon)}|| \epsilon + ||  T_{\alpha(\epsilon)}g - A^{\dagger} g || \rightarrow 0
\end{eqnarray*}
So,                                                                                 
$T_{\alpha(\epsilon)} g^{\epsilon}$ is close to 
$A^{\dagger} g$ provided $g^{\epsilon}$ is close to $g$.

The following three techniques are used for the regularization of ill-posed problems. 
\begin{itemize}
\item Truncated SVD:
$$ T_{\alpha}  = \sum_{\sigma_{k}\ge \alpha} \frac{1}{\sigma_{k}}
( \cdot  | \phi_k ) \psi_k \mbox{ .} $$
\item Tikhonov-Phillips regularization: 
$$ T_{\alpha} = (A^{\ast} A + \alpha I)^{-1} A^{\ast} \mbox{ ,} $$
which is associated to minimizing the quadratic form  
$$ || A f -  g^{\epsilon} ||^{2} + \alpha ||f||^{2} \mbox{ .}$$
More generally, it is associated to minimizing expressions of the form
$$ || A f -  g^{\epsilon} ||^{2} + \alpha q(f-f_0) \mbox{ ,}$$ where $q$ is some penalty function designed to keep $f$ close to a prior $f_0$.
\item Early stop of an iterative method: 
Assume that 
$$ f^{k+1} = B_{k} f^{k} + C_{k} g^{\epsilon} $$
is an iterative method to solve $A f = g $ with 
$B_k$ and $C_k$ bounded and  $\lim_{k\rightarrow\infty} f^{k} = A^{\dagger} g$.
For $\alpha > 0$, let the value of $k(\alpha)$ such that 
$k(\alpha )\rightarrow \infty$ when $\alpha\rightarrow 0$.
Then, under suitable conditions on $B_k$ and $C_k$, we have that $T_{\alpha} f \bydef f^{k(\alpha)}$ is a regularization of the problem.
\end{itemize}

In the case of infinite dimensional problems for compact operator equations of the form 
$A f = g$, it is natural to analyze how ill posed the problem is by looking at the rate of
convergence of the singular values of $A$ to zero. Problems for which the rate of decay of 
$\sigma_n$ is not faster than a polynomial are considered tractable. Problems for which
the decay is exponential or faster are considered severely ill-posed.
More details on regularization theory could be found for example in  \cite{EHN96, EKN89, ES00, TA77}. 

\subsection{Landweber - Kacmarz} \label{landwebb1}

We first consider Landweber's iteration to solve a nonlinear problem of the form 
$F(\gamma)=g$ where $F: D(F) \subset \H  \rightarrow \K$ is
a Fr\'{e}chet differentiable map between the Hilbert spaces $\H$ and $\K$. The iteration is defined by
\begin{equation}
\label{lbiteration}
\gamma_{k+1} = \gamma_{k} - F'(\gamma_k)^{\ast} (F(\gamma_k) - g)  \mbox{ ,} 
\end{equation}
where $\gamma_0 \in {D}(F)$

Under suitable conditions, the method converges \cite{BL05}. 
Furthermore, this iteration is known to generate a regularization method for the inverse
problem if we apply the early stopping technique mentioned above.
Landweber's iteration has been the subject of intense study both for theoretical as well as
for practical applications. See, for example, \cite{EHN96,ES00,HNS95}. 

In the sequel we describe the Landweber-Kaczmarz method for the doping profile
identification problem of Section~\ref{ssec:num-upm}. The notation and definitions of
the different operators follows that of Section~\ref{sec:6}.

\begin{itemize}
\item Parameter space: \ $\H \bydef  L^2(\Omega)$;
\item Input (fixed):  $ U_j \in H^{3/2}(\partial\Omega_D) \mbox{, }$ with 
$U_j |_{\Gamma_1} = 0$, \ $1 \le j \le N$; 
\item Output (data): \ $Y = \big\{ \Lambda_\gamma(U_j) \big\}_{j=1}^N
\in \K \bydef  \big[L^2(\Gamma_1)\big]^N $;
\item Parameter to output map: \ 
$F: \begin{array}[t]{rcl}
      {D} (F) \subset \H & \to & \K \\
      \gamma(x) & \mapsto & \big\{ \Lambda_\gamma(U_j) \big\}_{j=1}^N
    \end{array}$
\end{itemize}
where the domain of definition of the operator $F$ is
$$ D(F) \bydef \{ \gamma \in L^2(\Omega) ; \, \gamma_+ \ge \gamma(x)
           \ge \gamma_- > 0, \mbox{ a.e. in } \Omega \} \mbox{ .}$$
Here $\gamma_-$ and $\gamma_+$ are appropriate positive constants.
We shall denote the noisy data by $Y^\delta$ and assume that the data error
is bounded by $\| Y - Y^\delta \| \le \delta$.
Thus, we are able to represent the inverse doping problem in the form of finding 
\begin{equation} \label{eq:ip-dd}
F(\gamma) \, = \, Y^\delta .
\end{equation}  
It can be shown that \cite{BEMP01}
if we let the voltages $\{ U_j \}_{j=1}^N$ be chosen in the neighborhood of
$U \equiv 0$, then, the parameter-to-output map $F$ defined above is well-defined
and Fr\'echet differentiable on $D(F)$. Thus, 
the {\em Landweber iteration} \cite{DES98, EHN96, ES00, HNS95} becomes
$$ \gamma^\delta_{k+1} \ = \ \gamma_k^\delta - F'(\gamma_k^\delta)^*
   \big( F(\gamma_k^\delta) - Y^\delta \big)\, . $$
%

One possible variation of the Landweber iteration consists in coupling it
with the Kaczmarz strategy of considering an inner iteration where, at 
each step, one takes into account one component of the measurement vector only.
A detailed analysis of the Landweber-Kaczmarz method can be found in
\cite{KS02}. 

In the specific example mentioned above for Equation~(\ref{eq:ip-dd}), we
consider the components of the
parameter to output map: $F = \{ \mathcal F_j\}_{j=1}^N$, where
\begin{equation*}
\mathcal F_j: L^2(\Omega) \supset D(F) \ni \gamma \mapsto
\Lambda_\gamma(U_j) \in L^2(\Gamma_1) \, .
\end{equation*}
Now, setting $Y_j^\delta \bydef \mathcal F_j(\gamma^\delta)$,
$1 \le j \le N$, the Landweber-Kaczmarz iteration can be written in the form
\begin{equation} \label{eq:land-kacz}
  \gamma^\delta_{k+1} \, = \, \gamma_k^\delta
                          - \mathcal F_k'(\gamma_k^\delta)^*
  \big( \mathcal F_k (\gamma_k^\delta) - Y^\delta_k \big)\, ,
\end{equation}
for $k = 1, 2, \dots$, where we adopted the notation
$$  \mathcal F_k \bydef \mathcal F_j, \ \ Y^\delta_k \bydef Y^\delta_j,
    \ \ {\rm whenever} \ \ k = i \, N + j, \ \ {\rm and}\ \
    \left\{ \begin{array}{l}
    i = 0, 1, \dots \\ j = 1, \dots, N
    \end{array} \right. . $$

Notice that each step of the Landweber-Kaczmarz method consists of one
Landweber iterative step with respect to the $j$-th component of the
residual in (\ref{eq:ip-dd}). These Landweber steps are performed in a
cyclic way, using the components of the residual $\mathcal F_j(\gamma)
- Y^\delta_j$, $1 \le j \le N$, one at a time.

In the next section we shall describe how level set methods can be applied to tackle inverse problems.
A comparison between the Landweber-Kaczmarz method and a level set approach to 
the  doping profile identification problem was developed in \cite{LMZ05}. The preliminary
conclusion obtained therein is that in general the level set method performed better than
its Landweber-Kaczmarz counterpart.

\subsection{Level Set Methods in Inverse Problems} \label{levset}

The level set methodology has established itself as a promising alternative for the solution of several inverse problems that involve boundaries or obstacles.
The original formulation of level sets, as applied to curve and surface motion, is due to Osher and Sethian \cite{OS88}. The use 
of such methods in obstacle inverse problems is due to \cite{San95}. Burger \cite{Bur01} presented a rigorous mathematical
treatment for level set methods in inverse problems. See also \cite{LS03,FSL05} for a constrained optimization treatment of the 
method. In what follows we focus on the application of level set methods in inverse problems. 

Let $\Omega \subset \mathbb{R}^n $ be a given set and $F: \H \rightarrow \K $ a Fr\'echet differentiable operator.
The problem consists of finding $D \subset \mathrm{int}({\Omega}) $ in the equation 
\begin{equation} \label{levseteq1}
F(u) = g  \mbox{,}
\end{equation}
where 
$$ u = \left\{ 
\begin{array}{l}
u_{\mathrm{int}}, x\in D \\
u_{\mathrm{ext}}, x\in \Omega \setminus D
\end{array}
\right. $$
We now consider the boundary of the region $D$ in $\Omega$ as described by 
$\partial D = \left\{ x \in \Omega | \phi(x) = 0 \right\} \mbox{ .}$ The function $\phi$
shall be referred as the {\em level set function}. The level set function
evolves according to a parameter $t$ in such a way that
$$\partial D_{t}  \bydef  \left\{ x \in \Omega | \phi_{t} (x) = 0 \right\} \longrightarrow D   \mbox{ as } t \rightarrow \infty  \mbox{ .} $$
There are several possible dynamics for the evolution of $\phi_{t} $ with $t$. See, \cite{BL05} for
discussion and motivation, as well as \cite{Bur01,LS03}.  
One possibility is to use the dynamics introduced in
\cite{LS03,FSL05}. According to this approach, one represents  zero level set by an $H^1$-function $\phi: \Omega \to
\mathbb R$, in such a way that $\phi(x) > 0$ if $\gamma(x)=u_\mathrm{ext}$ and
$\phi(x) < 0$ if $\gamma(x)=u_\mathrm{int}$.
Starting from some initial guess $\phi_0 \in H^1(\Omega)$, one
solves the Hamilton-Jacobi equation
\begin{equation} \label{eq:ls-hj}
\frac{\partial \phi}{\partial t} + V \nabla \phi = 0
\end{equation}
where $V = v \frac{-\nabla\phi}{\ |\nabla\phi|^2}$ and the {\em velocity}
$v$ solves
\begin{equation} \label{eq:ls-vel}
\left\{ \!\! 
\begin{array}{l}
  \alpha (\Delta - I) v = \frac{\delta(\phi(t))}{|\nabla \phi(t)|}
  \left[F'(\chi(t))^*(F(\chi(t)) - Y^\delta) - \alpha\nabla \!\!\cdot\!\!
  \left( \frac{\nabla P(\phi)}{|\nabla P (\phi)|} \right) \right]
  ,\, {\rm in} \ \Omega \\
   \frac{\partial v}{\partial\nu} = 0 \ ,\, {\rm on}\ \partial\Omega
\end{array} \right.
\end{equation}
Here, $\alpha>0$ is a regularization parameter and $\chi = \chi(x,t)$ is
the projection of the level set function $\phi(x,t)$ defined by:
\begin{equation*}
\chi(x,t) = P(\phi(x,t)) \bydef
            \left\{ \begin{array}{ll}
              u_{\mathrm{ext}}, & {\rm if}\ \phi(x,t) > 0 \\
              u_{\mathrm{int}}, & {\rm if}\ \phi(x,t) < 0
            \end{array} \right. .
\end{equation*}

The above dynamics leads, in the case of the problem under consideration, to the following 

\noindent {\bf Algorithm:}
                                                                                
\begin{tt} \begin{enumerate}
\item[{\bf 1.}] \ Evaluate the residual 
$r_k \bydef F(P(\phi_k)) - Y^\delta$; 
\item[{\bf 2.}] \ Evaluate\, $v_k \bydef F'(P (\phi_k))^*(r_k)$;
\item[{\bf 3.}] \ Evaluate $w_k \in H^1(\Omega)$, satisfying
$$ \begin{aligned}
\alpha (I - \Delta) w_k &= - P'(\phi_k) \, v_k +
   \alpha P'(\phi_k) \nabla \cdot
   \left( \frac{\nabla P(\phi_k)}{|\nabla P (\phi_k)|} \right)\,,
   {\tt in}\ \Omega;\\
\frac{\partial v_k}{\partial \nu}|_{\partial\Omega}&=0\;.
\end{aligned} $$
\item[{\bf 4.}] \ Update the level set function \ $\phi_{k+1} =
\phi_k + \frac{1}{\alpha} \; v_k$.
\end{enumerate} \end{tt}
In practical implementations, instead of $P$, we use a smooth version $P_\varepsilon$.

\section{Some Numerical Experiments} \label{sec:6}

In this section we apply numerical methods to solve inverse doping profile
problems related to the V-C map. In the first two subsections, we address
the linearized unipolar case (close to equilibrium). In Subsection~%
\ref{ssec:num-upm} pointwise measurements of the current density are
considered, and in Subsection~\ref{ssec:num-ucfm} current flow measurements
through a contact are used as data. In the last subsection we present some
numerical results for the linearized bipolar case (close to equilibrium).

\subsection{Stationary linearized unipolar model: pointwise measurements of
the current density} \label{ssec:num-upm}

In this specific model, due to the assumptions $p = 0$ and $Q = 0$, the
Poisson equation and the continuity equation for the electron density
decouple. Therefore, we have to identify $C = C(x)$ from measurements of
the current density $\mu_n e^{V_{\rm bi}} \hat u_\nu |_{\Gamma_1}$,
where $(V^0, \hat u)$ solve, for each applied voltage $U$, the system
$$
\left\{ \begin{array}{rcl@{\ }l}
   \lambda^2 \, \Delta V^0 & = & e^{V^0} - C(x) & {\rm in}\ \Omega \\
   V^0 & = & V_{\rm bi}(x) & {\rm on}\ \partial\Omega_D \\
   \nabla V^0 \cdot \nu & = & 0 & {\rm on}\ \partial\Omega_N \\
\end{array} \right.
\hskip0.4cm
\left\{ \begin{array}{rcl@{\ }l}
   {\rm div}\, (\mu_n e^{V^0} \nabla \hat u) & = & 0 &  {\rm in}\ \Omega \\
   \hat u & = & U(x) & {\rm on}\ \partial\Omega_D \\
   \nabla \hat u \cdot \nu & = & 0 & {\rm on}\ \partial\Omega_N \, .
\end{array} \right.
$$

Notice that we split the problem in two parts: First we define the function
$\gamma(x) \bydef \mu_n e^{V^0(x)}$, $x \in \Omega$, and solve the parameter
identification problem
\begin{equation} \label{eq:num-d2n}
\left\{ \begin{array}{rcll}
   {\rm div}\, (\gamma \nabla \hat{u}) & = & 0 &  {\rm in}\ \Omega \\
   \hat{u} & = & U(x) & {\rm on}\ \Omega_D \\
   \nabla \hat{u} \cdot \nu & = & 0 & {\rm on}\ \Omega_N \, ,
\end{array} \right.
\end{equation}
for $\gamma$ from measurements of $\gamma \hat u_\nu |_{\Gamma_1}$. The
second step consists in the determination of $C$ in
$$ C(x) \ = \ \mu_n^{-1} \gamma(x) - \lambda^2 \,
   \Delta (\ln \mu_n^{-1}\gamma(x)) \, ,\ x \in \Omega \, . $$
The evaluation of $C$ from $\gamma$ is a mildly ill-posed problem and can be explicitly performed
in a routine way. We shall focus on the problem of identifying the function
parameter $\gamma$ in (\ref{eq:num-d2n}).
Therefore, the inverse doping profile problem in the linearized unipolar
model for pointwise measurements of the current density reduces to the
identification of the parameter $\gamma$ in (\ref{eq:num-d2n}) from
measurements of the Dirichlet to Neumann (DtN) map
$$  \Lambda_\gamma : \begin{array}[t]{rcl}
    H^{3/2}(\partial\Omega_D) & \to & H^{1/2}(\Gamma_1) \, . \\
    U & \mapsto & \gamma\, \hat u_\nu |_{\Gamma_1}
    \end{array} $$

If we take into account the restrictions imposed by the practical experiments
described in Subsection~\ref{ssec:ip-vcmap}, it follows:

{\em i)} The voltage profiles $U \in H^{3/2}(\partial\Omega_D)$ must
satisfy $U |_{\Gamma_1} = 0$;

{\em ii)}  The identification of $\gamma$ has to be performed from a finite
number of measurements, i.e. from the data
\begin{equation} \label{eq:data-upm}
\big\{ (U_j, \Lambda_\gamma(U_j)) \big\}_{j=1}^N
   \in \big[ H^{3/2}(\partial\Omega_D) \times H^{1/2}(\Gamma_1) \big]^N .
\end{equation}

\begin{figure}[b]
\centerline{
\includegraphics[width=5.8cm]{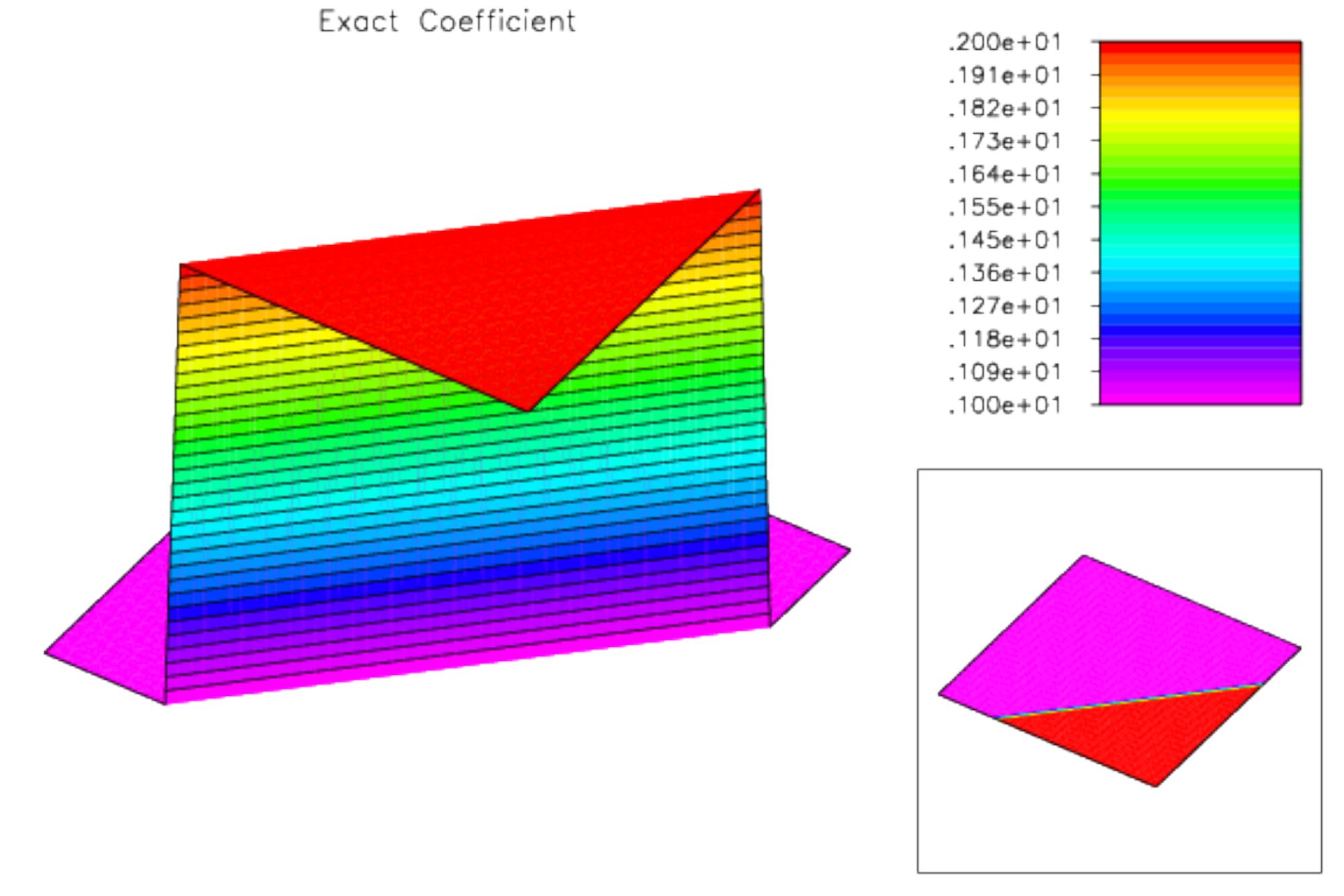}
\includegraphics[width=5.8cm]{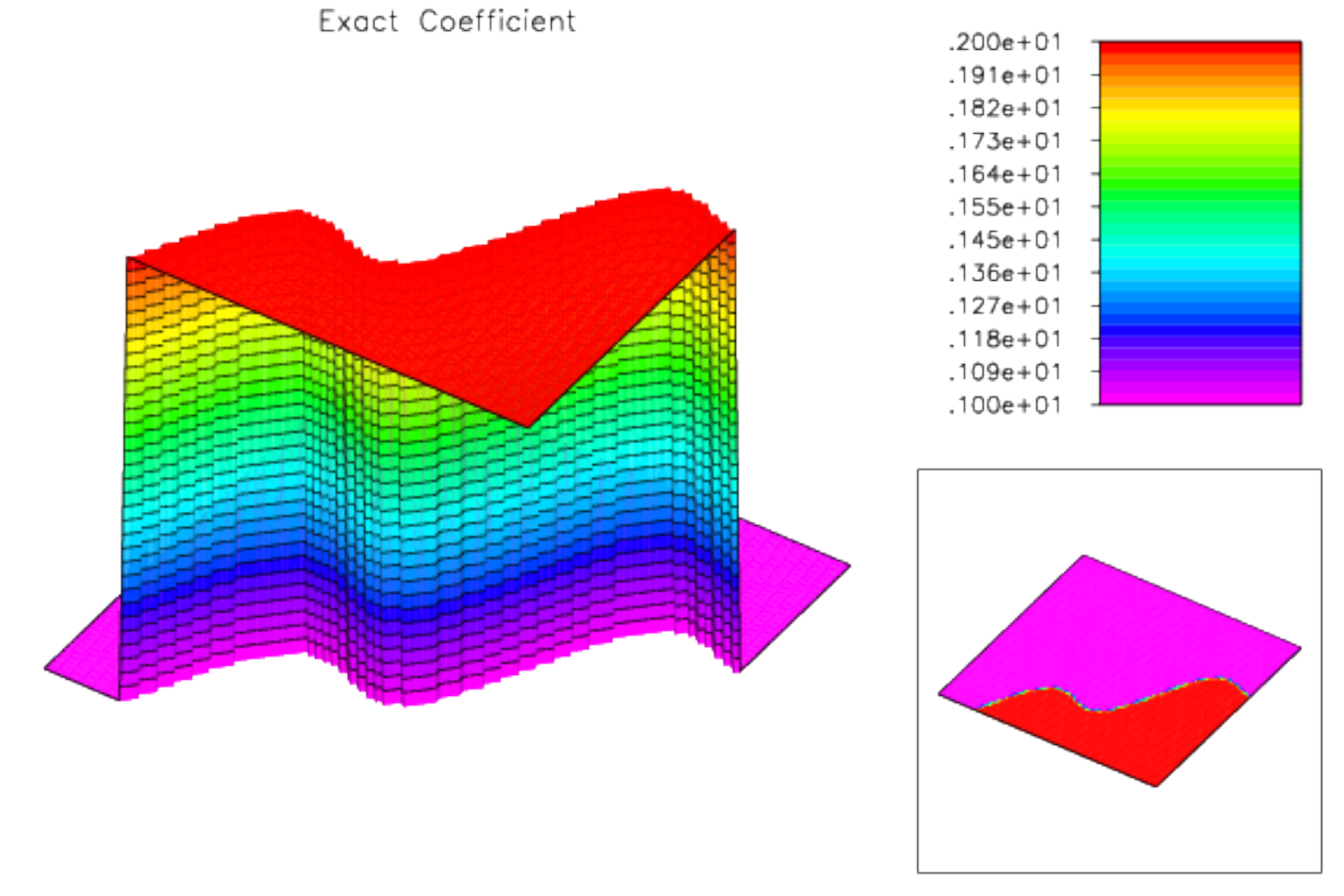} }
\centerline{\hfil (a) \hskip5cm  (b) \hfil}
\caption{\small Pictures (a) and (b) show the two different doping profiles
to be reconstructed in the numerical experiments.} \label{fig:exsol}
\end{figure}

For the concrete numerical tests presented in this paper, we apply an
iterative method of level set type to solve problem (\ref{eq:ip-abstract})
See \cite{LMZ05} for details.
The domain $\Omega \subset \mathbb R^2$ is the unit square, and the boundary
parts are defined as follows
$$ \Gamma_1 \ \bydef \  \{ (x,1) \, ;\ x \in (0,1) \} \, , \ \
   \Gamma_0 \ \bydef \  \{ (x,0) \, ;\ x \in (0,1) \} \, , $$
$$ \partial\Omega_N \ \bydef \ \{ (0,y) \, ;\ y \in (0,1) \} \cup 
   \{ (1,y) \, ;\ y \in (0,1) \} \, . $$
The fixed inputs $U_j$, are chosen to be piecewise constant functions
supported in $\Gamma_0$
$$  U_j(x) \ \bydef \ \left\{ \begin{array}{rl}
      1, & |x - x_j| \le h \\
      0, & {\rm else} \end{array} \right. $$
where the points $x_j$ are equally spaced in the interval $(0,1)$. The
doping profiles to be reconstructed are shown in Figure~\ref{fig:exsol}.
In these pictures, as well as in the forthcoming ones, $\Gamma_1$ is the
lower left edge and $\Gamma_0$ is the top right edge (the origin corresponds
to the upper right corner).

For the experiments concerning pointwise measurements of the current
density, we assume that only one measurement is available, i.e. $N = 1$
in (\ref{eq:data-upm}).

\begin{figure}[t]
\centerline{ \includegraphics[width=8.4cm]{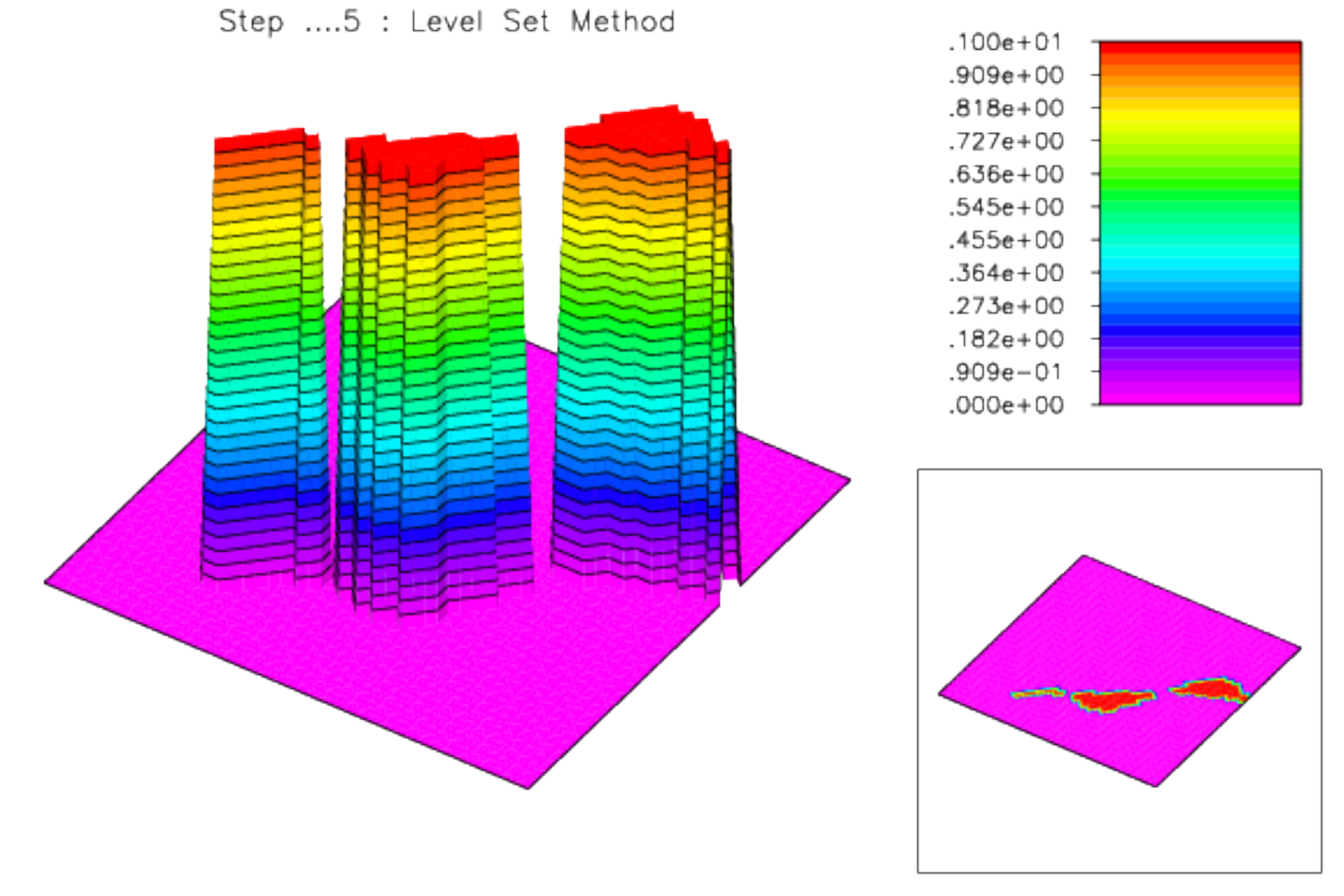} }
\centerline{ \includegraphics[width=8.4cm]{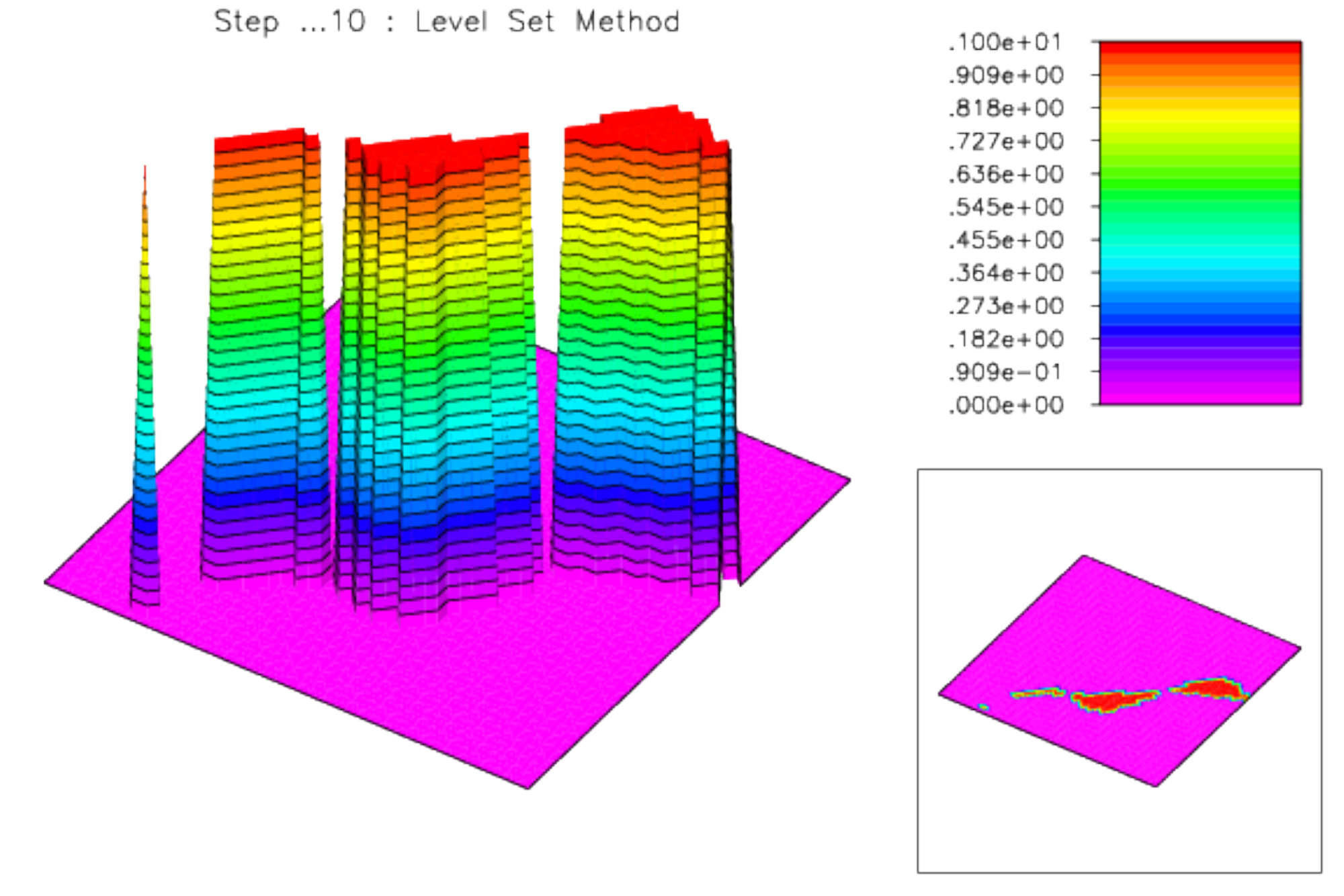} }
\centerline{ \includegraphics[width=8.4cm]{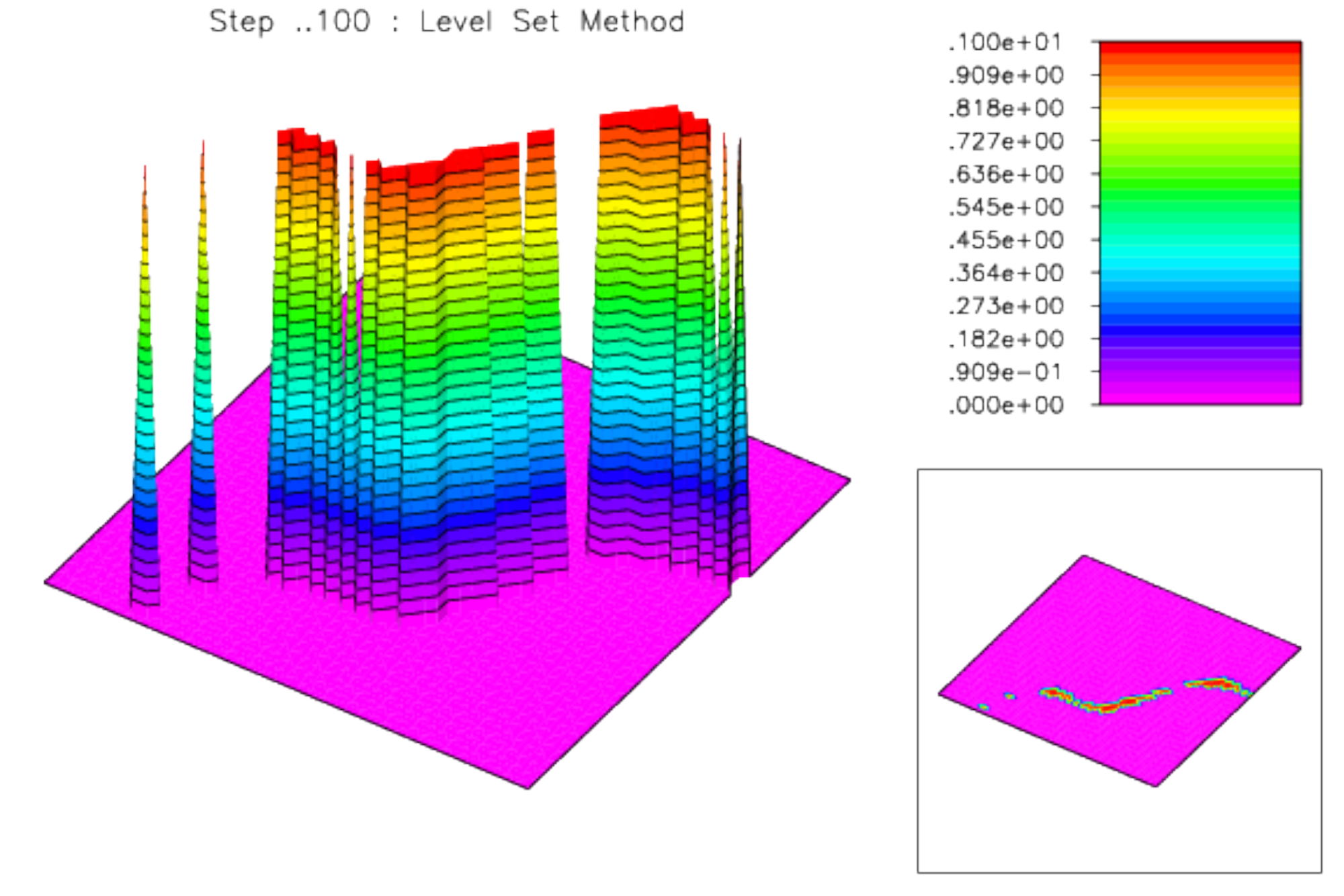} }
\caption{\small First experiment for the unipolar model with pointwise
measurements of the current density:
Reconstruction of the p-n junction in Figure~\ref{fig:exsol}~(b).
Evolution of the iteration error for exact data and one measurement of
the DtN map $\Lambda_\gamma$ (i.e. $N = 1$ in (\ref{eq:data-upm})).}
\label{fig:exp1-up-pm}
\end{figure}

The first numerical experiment is shown in Figure~\ref{fig:exp1-up-pm}.
Here exact data is used for the reconstruction of the p-n junction in
Figure~\ref{fig:exsol}~(b). The pictures correspond to plots of the
iteration error after 5, 10 and 100 steps of the level set method.

\begin{figure}[t]
\centerline{ \includegraphics[width=8.4cm]{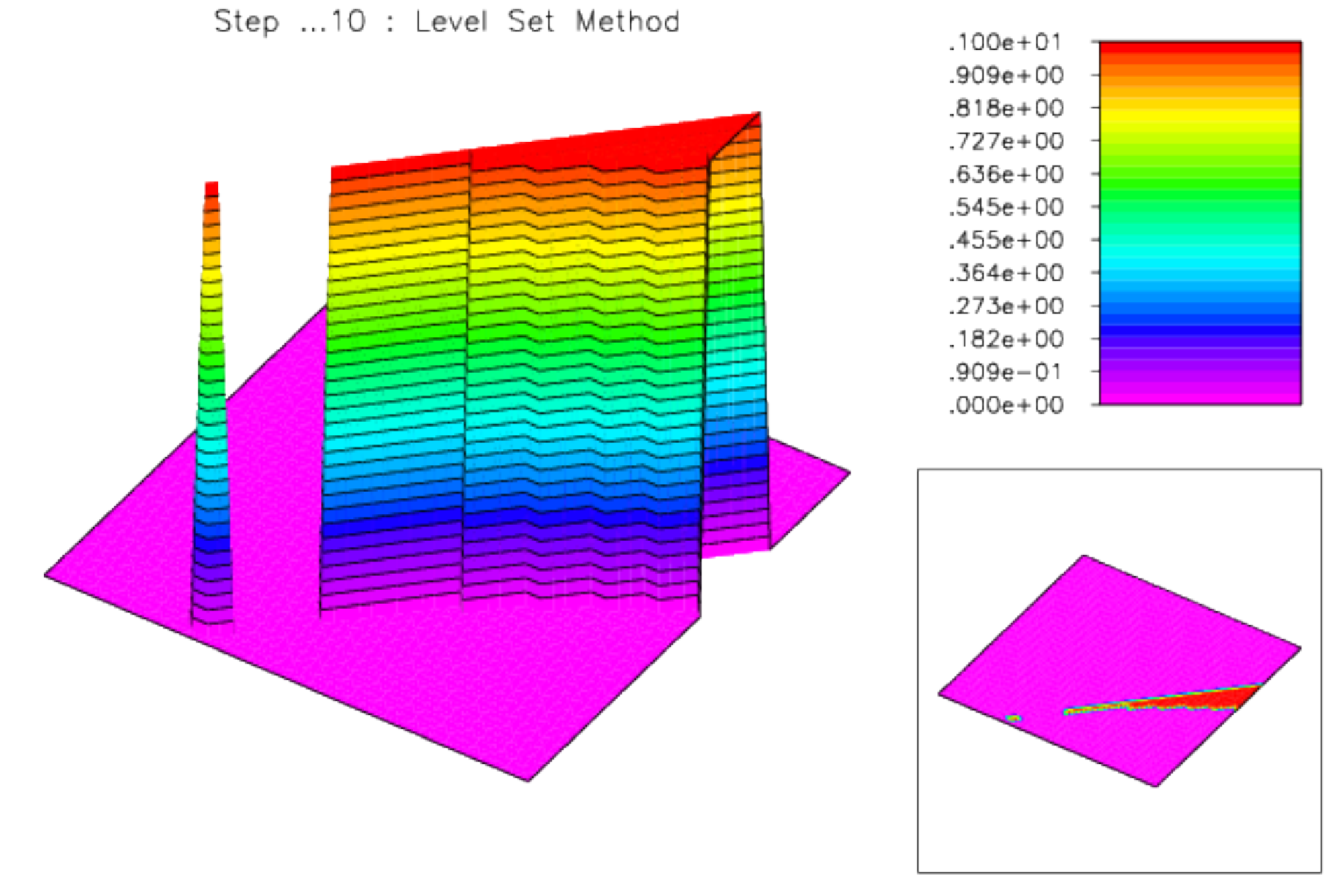} }
\centerline{ \includegraphics[width=8.4cm]{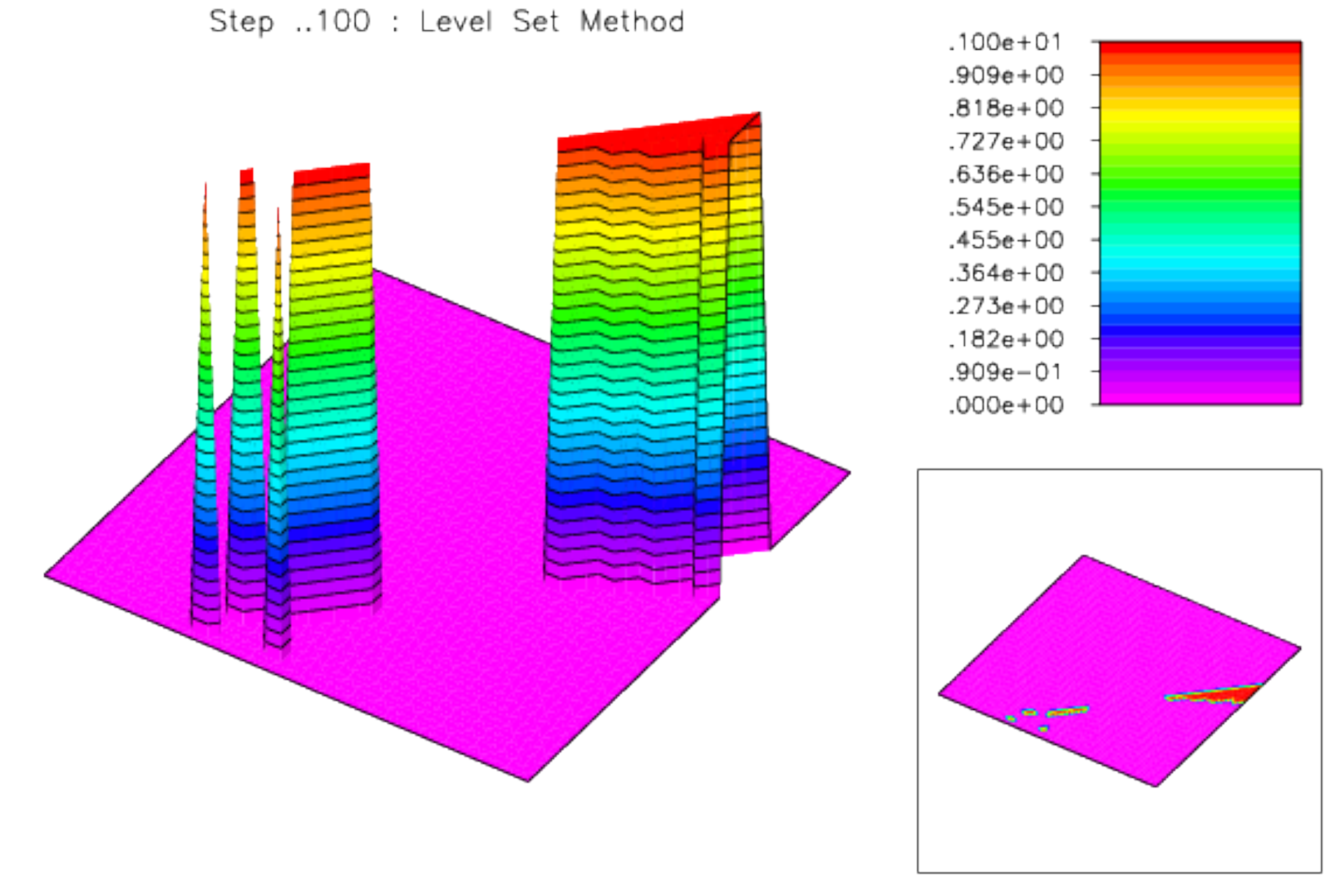} }
\centerline{ \includegraphics[width=8.4cm]{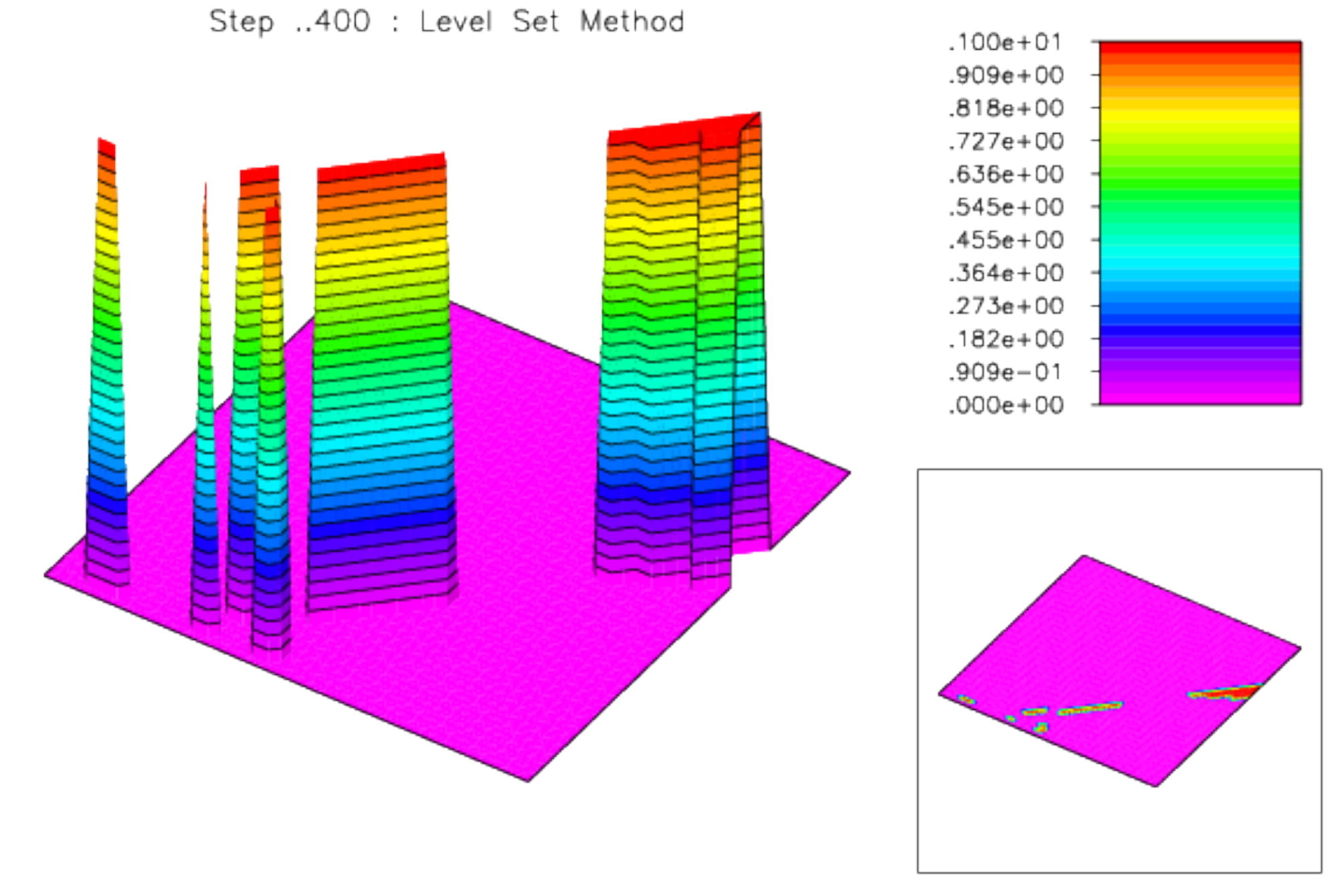} }
\caption{\small Second experiment for the unipolar model with pointwise
measurements of the current density:
Reconstruction of the p-n junction in Figure~\ref{fig:exsol}~(a).
Only one measurement of the DtN map $\Lambda_\gamma$ is available
(i.e. $N = 1$ in (\ref{eq:data-upm})).
Evolution of the iteration error for the level set method and data
contaminated with 10\% random noise.} \label{fig:exp2-up-pm}
\end{figure}

The second experiment (see Figure~\ref{fig:exp2-up-pm}) concerns the
reconstruction of the p-n junction in Figure~\ref{fig:exsol}~(a). In
this experiment the data is contaminated with 10\% random noise. The
pictures correspond to plots of the iteration error after 10, 100 and 400
steps of the level set method.

\subsection{Stationary linearized unipolar model: current flow measurements
through a contact} \label{ssec:num-ucfm}

In what follows we consider the same unipolar model as in Subsection~%
\ref{ssec:num-upm}. Again we shall focus on the identification problem
related to (\ref{eq:num-d2n}). However, the coefficient $\gamma$ has to
be identified from measurements of the current flow through the contact
$\Gamma_1$, i.e. from
$$ \int_{\Gamma_1} \gamma \hat u_\nu \, ds \, , $$
where $\hat u$ solve (\ref{eq:num-d2n}) for prescribed inputs $U \in
H^{3/2}(\partial\Omega_D)$.

An immediate remark is that the amount of available data is much larger
in the case of pointwise measurements of the current density than in the case
of current flow measurements through a contact.
Notice that the inverse doping profile problem in the linearized unipolar
model for measurements of the current flow through the contact $\Gamma_1$
reduces to the identification of the parameter $\gamma$ in (\ref{eq:num-d2n})
from measurements of the (averaged) DtN map
$$  \widetilde\Lambda_\gamma : \begin{array}[t]{rcl}
    H^{3/2}(\partial\Omega_D) & \to & \mathbb R \\
    U & \mapsto & \int_{\Gamma_1} \gamma\, \hat u_\nu \, ds
    \end{array} $$

As in the previous subsection, we take into account the restrictions
imposed by practical experiments, which lead to the following assumptions:

{\em i)} The voltage profile $U \in H^{3/2}(\partial\Omega_D)$ must
satisfy $U |_{\Gamma_1} = 0$;

{\em ii)}  The identification of $\gamma$ has to be performed from a finite
number of measurements, i.e. from the data
\begin{equation} \label{eq:data-cfm}
\big\{ (U_j, \widetilde\Lambda_\gamma(U_j)) \big\}_{j=1}^N 
       \in \big[ H^{3/2}(\partial\Omega_D) \times \mathbb R \big]^N .
\end{equation}

The subsequent numerical tests were performed using the same iterative
method of level set type as in the previous subsection. The domain
$\Omega \subset \mathbb R^2$ as well as the boundary parts $\Gamma_0$,
$\Gamma_1$ and $\partial\Omega_N$ are defined as before.

For the experiments concerning current flow measurements through the
contact $\Gamma_1$, we assume that several measurements are available,
i.e. $N > > 1$ in (\ref{eq:data-cfm}).

\begin{figure}[t]
\centerline{ \includegraphics[width=8.4cm]{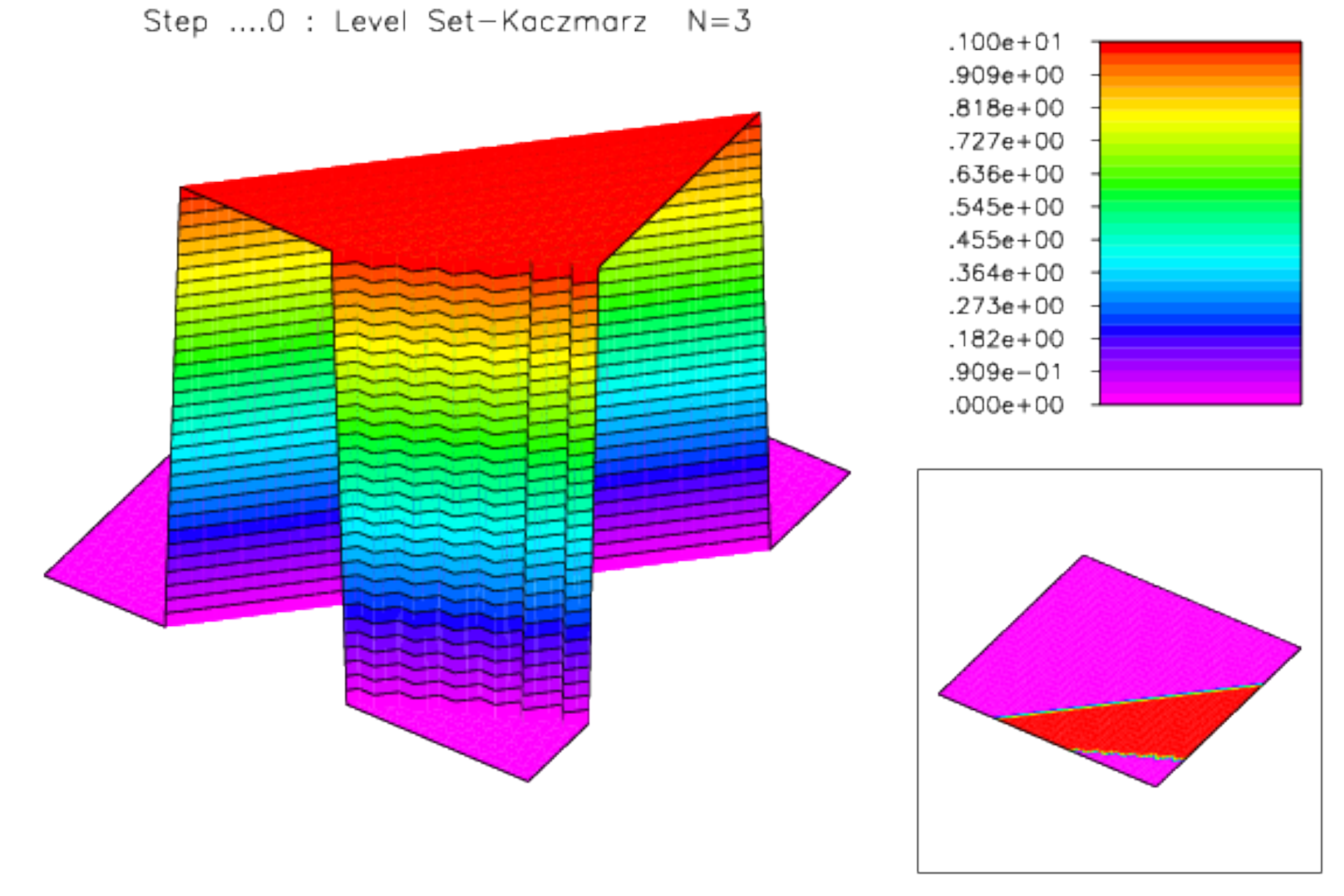} }
\centerline{ \includegraphics[width=8.4cm]{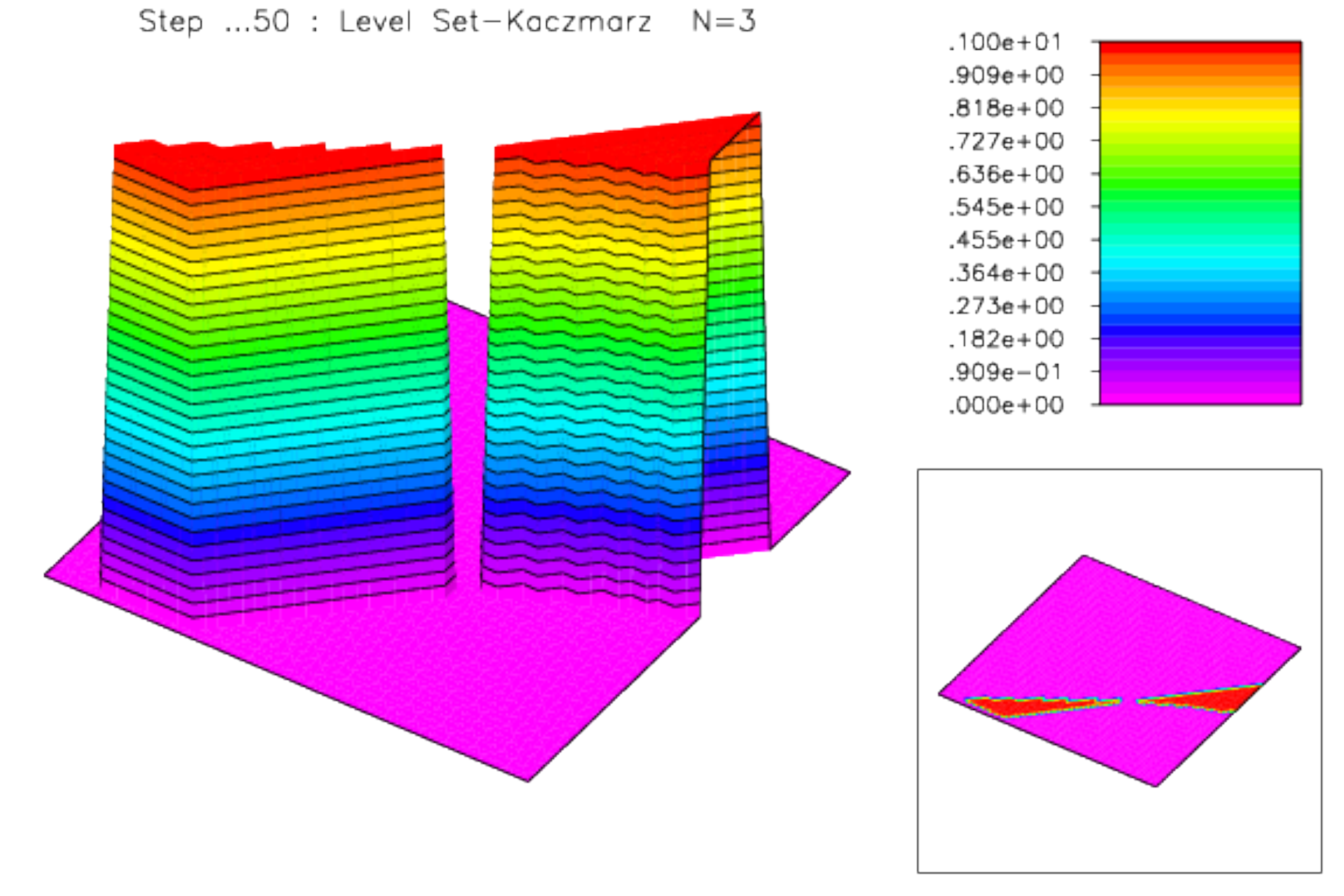} }
\centerline{ \includegraphics[width=8.4cm]{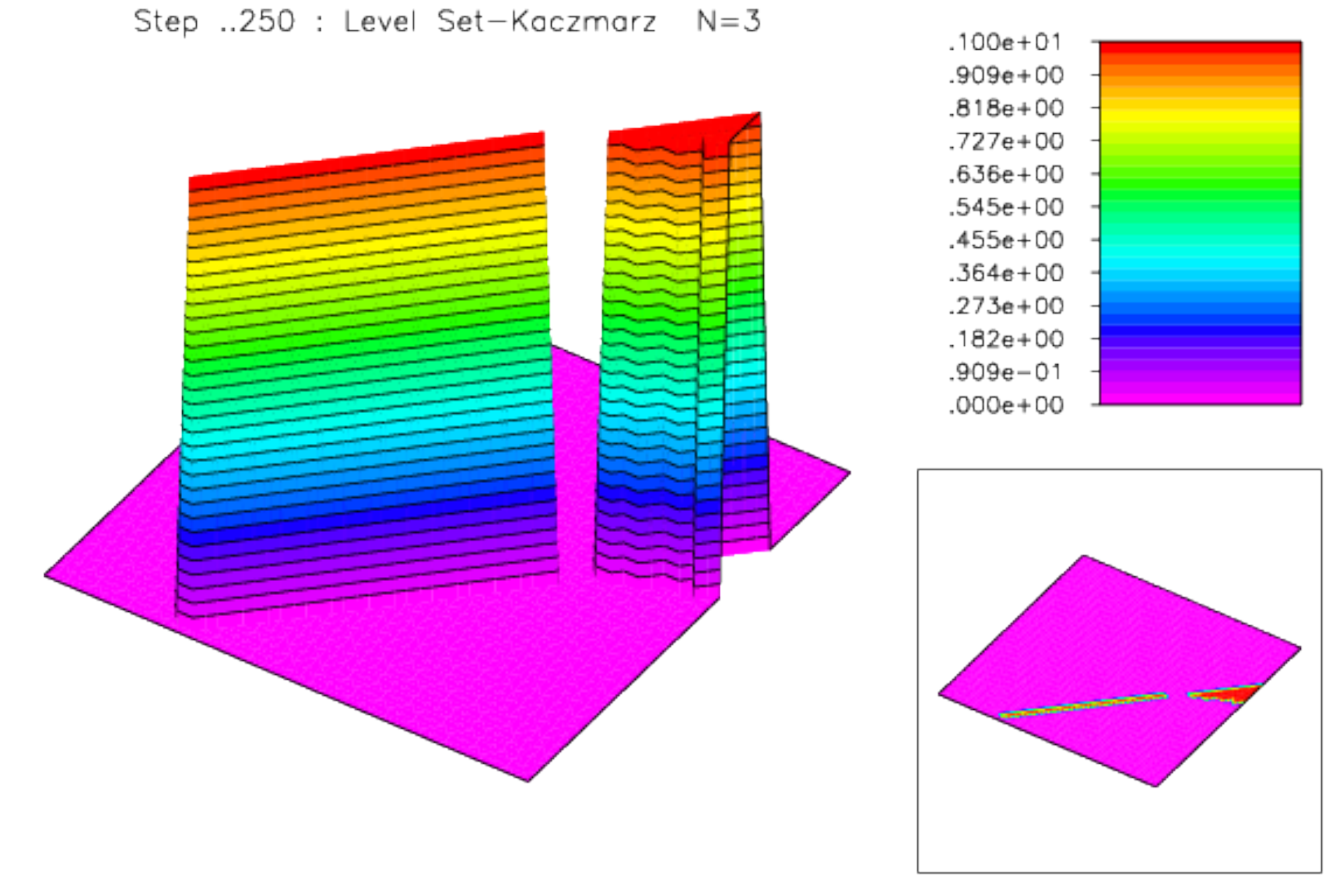} }
\caption{\small First experiment for the unipolar model with current flow
measurements through the contact $\Gamma_1$:
Reconstruction of the p-n junction in Figure~\ref{fig:exsol}~(a).
Three measurements of the DtN map $\widetilde\Lambda_\gamma$ are used in
the reconstruction (i.e. $N = 3$ in (\ref{eq:data-cfm})).
Evolution of the iteration error for exact data.} \label{fig:exp1-up-cfm}
\end{figure}

The first numerical experiment is shown in Figure~\ref{fig:exp1-up-cfm}.
Here exact data is used for the reconstruction of the p-n junction in
Figure~\ref{fig:exsol}~(a). The picture on the left hand side shows the
error for the initial guess of the iterative method.%
\footnote{In all numerical experiments presented in this paper we used
the same initial guess for the iterative methods. We observed that the
choice of the initial guess does not significantly influence the overall
performance of the iterative method.}
The other two pictures correspond to plots of the iteration error after
50 and 250 steps of the level set method respectively.

\begin{figure}[t]
\centerline{ \includegraphics[width=8.4cm]{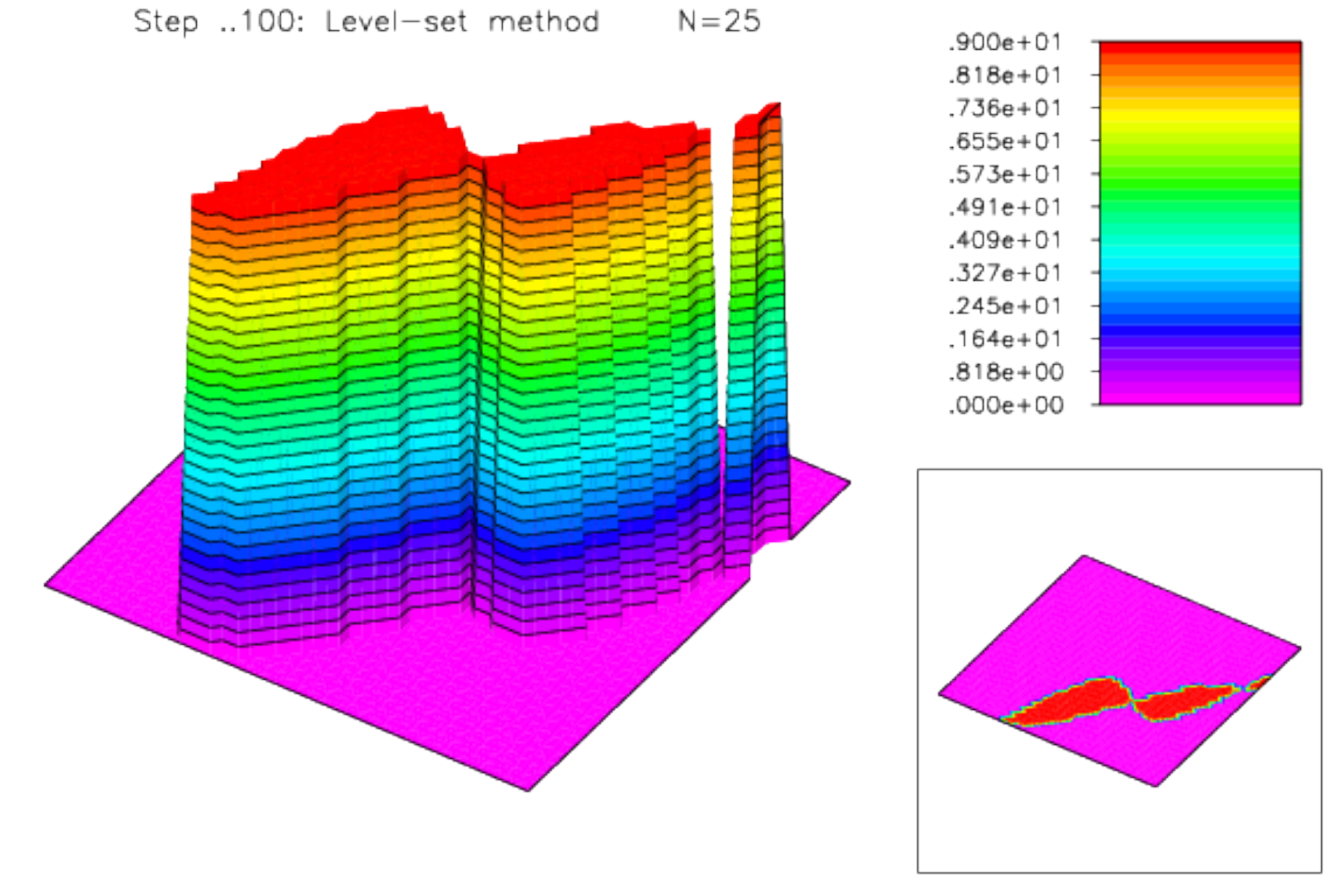} }
\centerline{ \includegraphics[width=8.4cm]{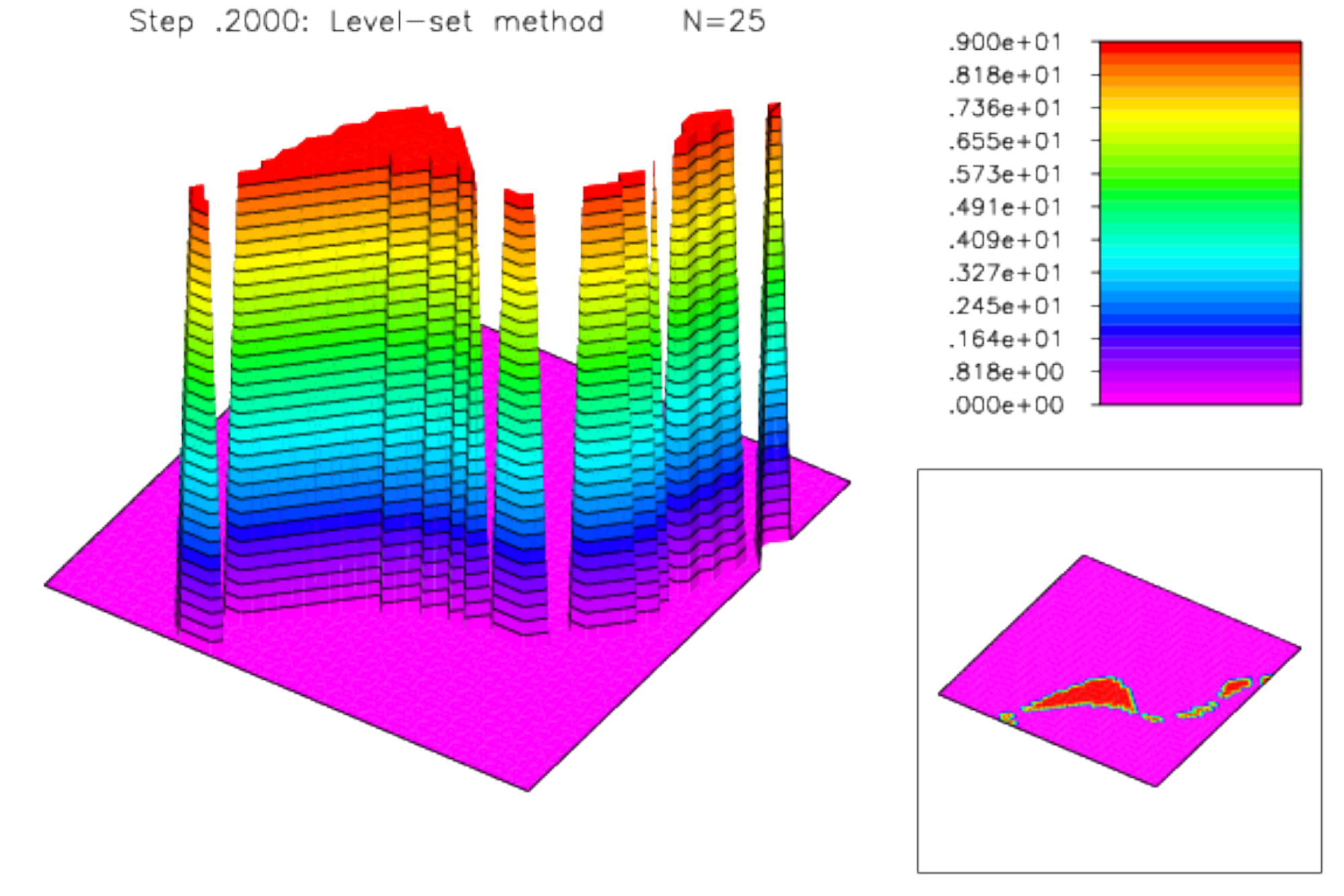} }
\centerline{ \includegraphics[width=8.4cm]{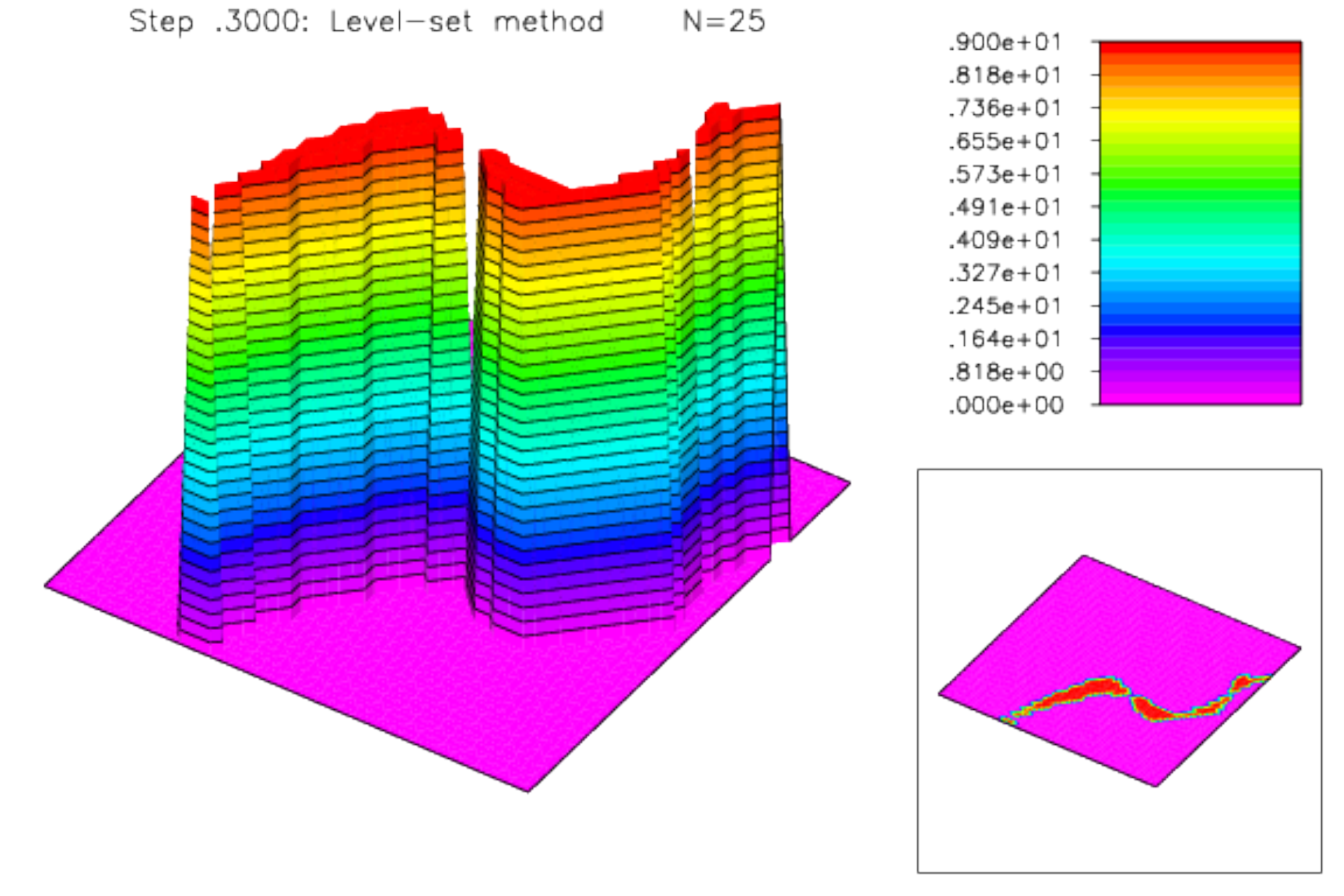} }
\caption{\small Second experiment for the unipolar model with current flow
measurements through the contact $\Gamma_1$:
Reconstruction of the p-n junction in Figure~\ref{fig:exsol}~(b).
The data consists of twenty five measurements of the DtN map
$\widetilde\Lambda_\gamma$ (i.e. $N = 25$ in (\ref{eq:data-cfm})).
Plots of the iteration error after 100, 2000 and 3000 steps. Data
with 1\% random noise.} \label{fig:exp2-up-cfm}
\end{figure}

The second experiment (see Figure~\ref{fig:exp2-up-cfm}) concerns
the reconstruction of the p-n junction in Figure~\ref{fig:exsol}~(b).
The available data is contaminated with 1\% random noise. The
pictures correspond to plots of the iteration error after 100, 2000
and 3000 steps of the level set method.

\subsection{Stationary linearized bipolar model: pointwise measurements of
the current density} \label{ssec:num-bpm}

In the sequel we consider the bipolar model introduced in Subsection~%
\ref{ssec:bipol}. As in the unipolar model, it follows from the assumption
$Q = 0$ that the Poisson equation (\ref{eq:equil-caseA}) and the continuity
equations (\ref{eq:bipol-statA}), (\ref{eq:bipol-statB}) decouple.
The inverse doping profile problem corresponds to the identification of
$C = C(x)$ from pointwise measurements of the total current density
$J = J_n + J_p$, namely
$$ ( \mu_n e^{V_{\rm bi}} \hat{u}_\nu - \mu_p e^{-V_{\rm bi}} \hat{v}_\nu )
|_{\Gamma_1}  \mbox{ .}$$
Compare with the Gateaux derivative of the V-C map $\Sigma_C$ at the point
$U=0$ in (\ref{eq:def-sigma-prime-C}).
Here $(V^0, \hat u, \hat v)$ solve, for each applied voltage $U$, the system
(\ref{eq:equil-case}), (\ref{eq:bipol-stat}) (with $h$ substituted by $U$).

As in the unipolar case, we can split the inverse problem in two parts:
First we define the function $\gamma(x) \bydef e^{V^0(x)}$, $x \in \Omega$,
and solve the parameter identification problem
\begin{equation} \label{eq:num-d2nB}
\left\{ \begin{array}{rcl@{\ }l}
   {\rm div}\, (\mu_n \gamma \nabla \hat u) & = & 0 &  {\rm in}\ \Omega \\
   \hat u & = & - U(x) & {\rm on}\ \partial\Omega_D \\
   \nabla \hat u \cdot \nu & = & 0 & {\rm on}\ \partial\Omega_N
\end{array} \right.
\hskip0.2cm
\left\{ \begin{array}{rcl@{\ }l}
   {\rm div}\, (\mu_p \gamma^{-1} \nabla \hat v) & = & 0 &  {\rm in}\ \Omega \\
   \hat v & = & U(x) & {\rm on}\ \partial\Omega_D \\
   \nabla \hat v \cdot \nu & = & 0 & {\rm on}\ \partial\Omega_N
\end{array} \right.
\end{equation}
for $\gamma$, from measurements of 
$( \mu_n \gamma \hat{u}_\nu - \mu_p \gamma^{-1} \hat{v}_\nu )
|_{\Gamma_1}$. The second step consists in the determination of $C$ in
$$ C(x) \ = \ \gamma(x) - \gamma^{-1}(x) - \lambda^2 \,
   \Delta (\ln \gamma(x)) \, ,\ x \in \Omega \, . $$

Analogous to the unipolar case, the evaluation of $C$ from $\gamma$ can be
performed in a stable way. Therefore, we shall focus on the problem of identifying
the function parameter $\gamma$ in (\ref{eq:num-d2nB}).
Notice that the inverse doping profile problem in the linearized bipolar
model for pointwise measurements of the current density reduces to the
identification of the parameter $\gamma$ in (\ref{eq:num-d2nB}) from
measurements of the Dirichlet to Neumann (DtN) map
$$  \Phi_\gamma : \begin{array}[t]{rcl}
    H^{3/2}(\partial\Omega_D) & \to & H^{1/2}(\Gamma_1) \, . \\
    U & \mapsto & ( \mu_n \gamma \hat{u}_\nu - \mu_p \gamma^{-1} \hat{v}_\nu )
                  |_{\Gamma_1}
    \end{array} $$

As before we take into account the restrictions imposed by the practical
experiments, from what follows:

{\em i)} The voltage profiles $U \in H^{3/2}(\partial\Omega_D)$ must
satisfy $U |_{\Gamma_1} = 0$;

{\em ii)}  The identification of $\gamma$ has to be performed from a finite
number of measurements, i.e. from the data
\begin{equation} \label{eq:data-bp-pm}
\big\{ (U_j, \Phi_\gamma(U_j)) \big\}_{j=1}^N
   \in \big[ H^{3/2}(\partial\Omega_D) \times H^{1/2}(\Gamma_1) \big]^N .
\end{equation}

\begin{figure}[t]
\centerline{ \includegraphics[width=8.4cm]{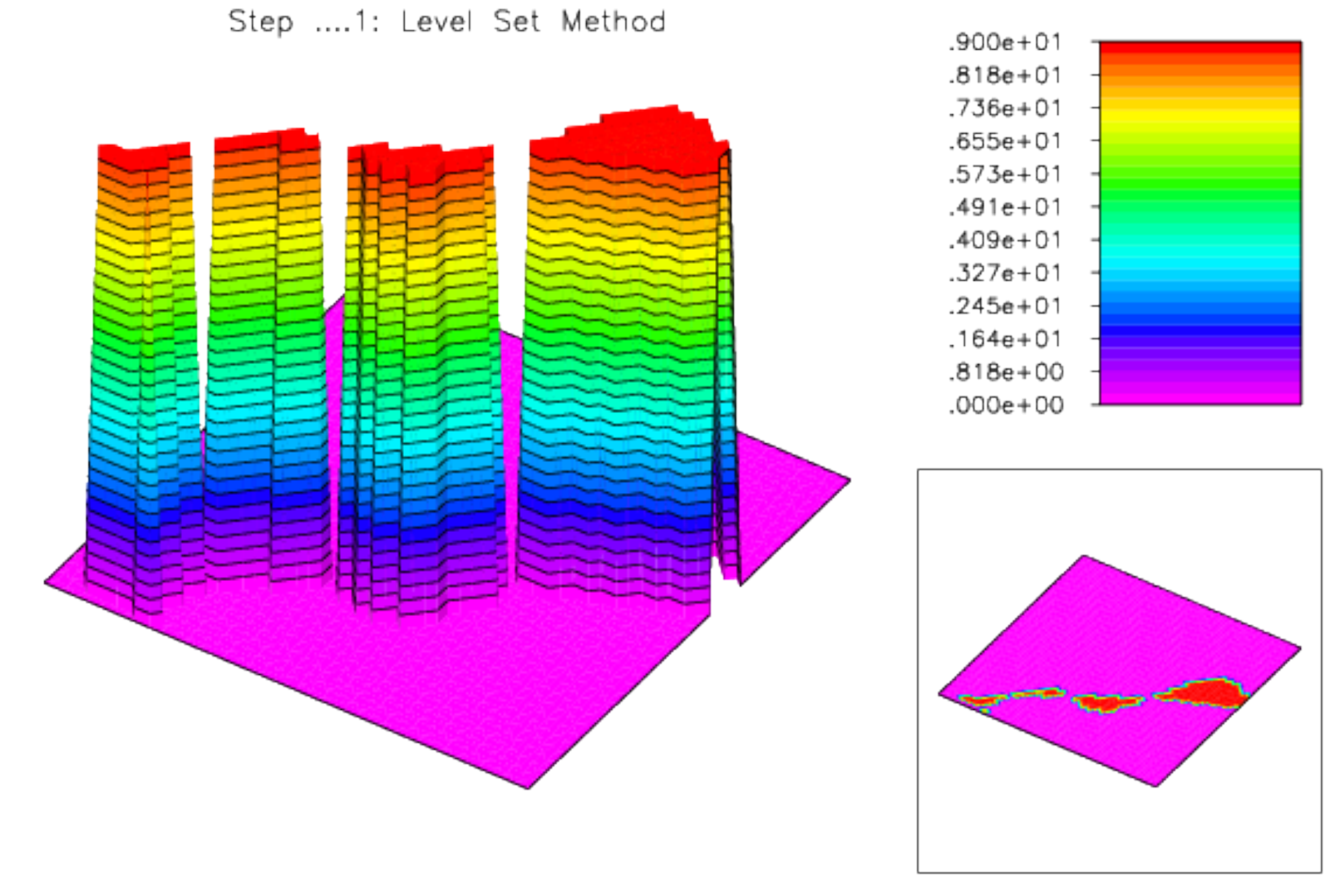} }
\centerline{ \includegraphics[width=8.4cm]{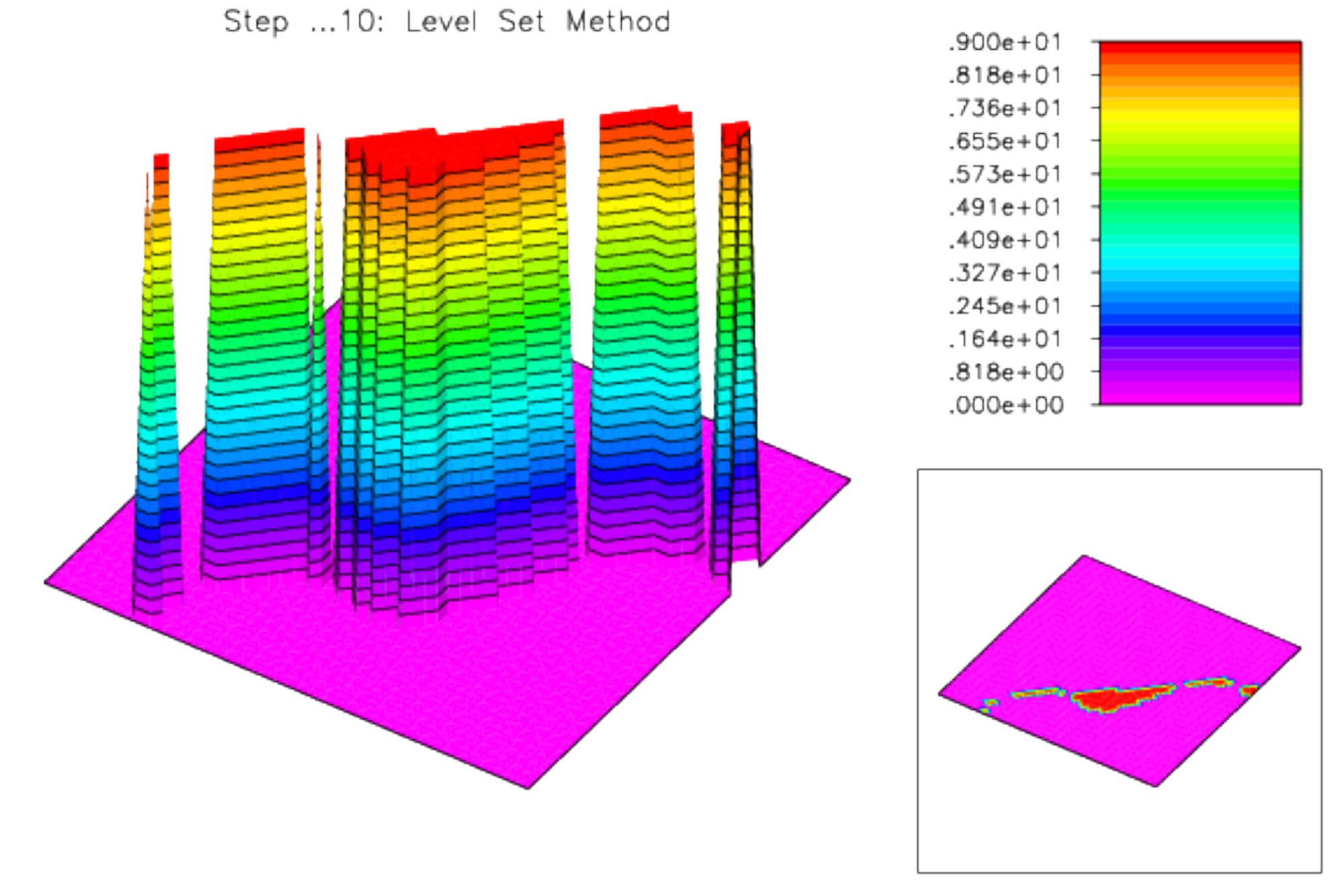} }
\centerline{ \includegraphics[width=8.4cm]{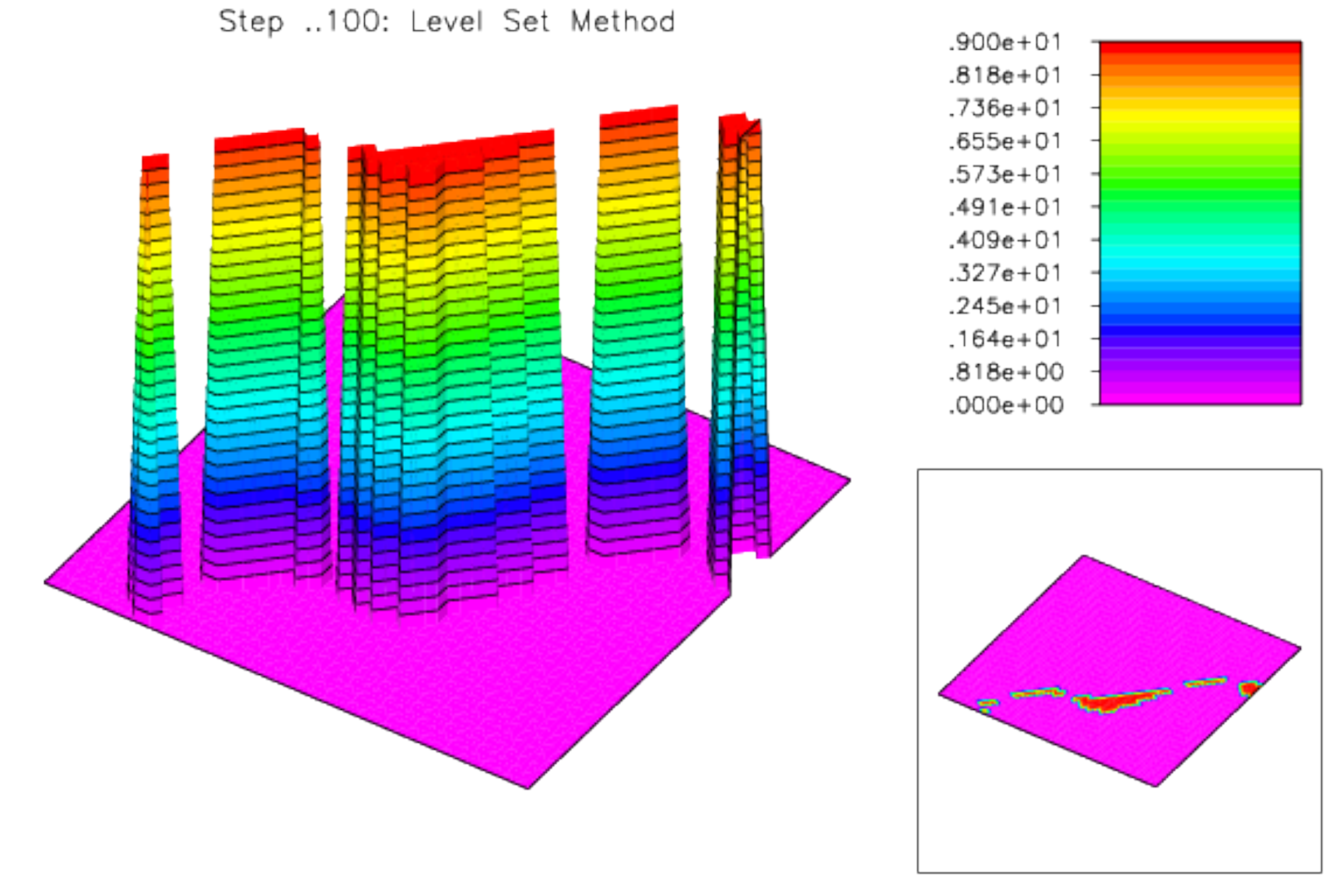} }
\caption{\small Experiment for the bipolar model with pointwise measurements
of the current density:
Reconstruction of the p-n junction in Figure~\ref{fig:exsol}~(b).
Evolution of the iteration error for exact data and one measurement of
the DtN map $\Phi_\gamma$ (i.e. $N = 1$ in (\ref{eq:data-bp-pm})).}
\label{fig:exp-bp-pm}
\end{figure}

In Figure~\ref{fig:exp-bp-pm} we present a numerical experiment for the
bipolar model with pointwise measurements of the current density. Here
exact data is used for the reconstruction of the p-n junction in Figure~%
\ref{fig:exsol}~(b). The pictures show plots of the iteration error after 1,
10 and 100 steps of the level set method respectively.

\subsection{Remarks and conclusions}

The best numerical results are obtained for the experiments concerning
the linearized unipolar case with pointwise measurements of the current
density.
In this model, a single measurement of the DtN map $\Lambda_\gamma$,
i.e. $N = 1$ in (\ref{eq:data-upm}), contains enough information about
the structure of the doping profile and suffices to obtain a very precise
reconstruction of the p-n junction. This is the case even for highly
oscillating p-n junctions as shown in  Figure~\ref{fig:exp1-up-pm} and also in
the presence of noise (see Figure~\ref{fig:exp2-up-pm}). We observed that
the iteration is extremely robust with respect to the choice of the
initial guess and also with respect to high levels of noise. 

Concerning the linearized unipolar model with current flow measurements
through the contact $\Gamma_1$, our experiments show that the (averaged)
DtN map $\widetilde\Lambda_\gamma$, furnishes much less information about
the solution structure than the map $\Lambda_\gamma$.
Depending on the complexity of the p-n junction, more measurements of
$\widetilde\Lambda_\gamma$ may be needed in order to obtain an acceptable
reconstruction.
The experiments show that a single measurement ($N = 1$ in
(\ref{eq:data-cfm})) is not enough to identify the doping profile in
Figure~\ref{fig:exsol}~(a). Moreover, although three measurements have
shown to be enough to reconstruct this p-n junction (see Figure~%
\ref{fig:exp1-up-cfm}), this is not the case for the p-n junction in
Figure~\ref{fig:exsol}~(b). For this second and more complex junction,
we first obtained more accurate reconstructions with $N = 19$ in
(\ref{eq:data-cfm}). The quality of the reconstruction obtained for
$N = 25$ is already very high (see Figure~\ref{fig:exp2-up-cfm}) and
does not qualitatively improve for larger values of $N$ (we experimented
up to $N = 49$).

It is worth noticing that the number of iterative steps required by the
level set algorithm to reach the stopping criteria for the inverse problem
related to the map $\widetilde\Lambda_\gamma$ is greater than that for the
operator $\Lambda_\gamma$. This is again explained by the fact that the
range of $\widetilde\Lambda_\gamma$ lies in $\mathbb R$, while the range
of $\Lambda_\gamma$ lies in $H^{1/2}(\Gamma_1)$.

Concerning the experiments for the linearized bipolar model with pointwise
measurements of the current density, the quality of the results is comparable
to those in Subsection~\ref{ssec:num-upm} and, as in that subsection, a single
measurement of the operator $\Phi_\gamma$ ($N = 1$ in (\ref{eq:data-bp-pm}))
suffices to precisely reconstruct the p-n junction. We observed, however,
that convergence of the iterative method is more sensitive to the choice of
the initial condition than in Subsection~\ref{ssec:num-upm}.

\section*{Appendix}

Properties of silicon at room temperature

\begin{table}
\begin{center} \begin{tabular}{cl}
\hline \\[-2.3ex]
{ Parameter} & { \hfil Typical value \hfil} \\[0.3ex]
\hline \\[-2.0ex]
$\epsilon_s$ & $11.9 \ \epsilon_0$ \\
$\mu_n$      & $\approx 1500 \ {\rm cm}^2 \ {\rm V}^{-1} \ {\rm s}^{-1}$ \\
$\mu_p$      & $\approx  450 \ {\rm cm}^2 \ {\rm V}^{-1} \ {\rm s}^{-1}$ \\
$C_n$    & $2.8 \times 10^{-31} \ {\rm cm}^6 {\rm / s}$ \\
$C_p$    & $9.9 \times 10^{-32} \ {\rm cm}^6 {\rm / s}$ \\
$\tau_n$ & $10^{-6}\, {\rm s}$ \\
$\tau_p$ & $10^{-5}\, {\rm s}$ \\
\hline
\end{tabular} \end{center}
\caption{Typical values of main the constants in the model.
\label{tab:typ-val}}
\end{table}

\noindent
Relevant physical constants:

\begin{itemize}
\item[] Permittivity of vacuum:
       $\epsilon_0 = 8.85 \times 10^{-14} {\rm As \, V}^{-1} \, {\rm cm}^{-1}$;
\item[] Elementary charge:
       $q = 1.6 \times 10^{-19} {\rm As}$.
\end{itemize}

\section*{Acknowledgment}

AL acknowledges support from the Brazilian National Research Council CNPq,
under project grants 305823/03-5 and 478099/04-5.
PAM acknowledges support from the Austrian National Science Foundation
FWF through his Wittgenstein Award 2000.
J.P.Z. acknowledges financial support from CNPq through grants 302161/2003-1 and
474085/2003-1.


\printindex

\end{document}